\documentclass[invmat]{svjour}
\usepackage{psnfss/times}

\usepackage{amsfonts}
\usepackage{amssymb}
\usepackage{amsmath}
\usepackage{theorem}

\newtheorem{theorem}{Theorem}[section]
\newtheorem{lemma}[theorem]{Lemma}

\newtheorem{proposition}[theorem]{Proposition}
\newtheorem{corollary}[theorem]{Corollary}
 {\theorembodyfont{\rmfamily}}

\def \RR {\mathbb R}

\def \EE {\mathbb E}

\def \eps {\varepsilon}
\def \vphi {\varphi}

\def \E {\mathcal E}

\begin{document}

\title{A Central Limit Theorem for Convex Sets}
\author{B. Klartag\thanks{The author is a Clay Research Fellow
and is also supported by NSF grant $\#DMS-0456590$.}} \institute{
Department of Mathematics, Princeton University, Princeton, NJ
08540, USA}
\date{}

\maketitle

\renewcommand{\theequation}{$\ast$}
\begin{abstract}
We show that there exists a sequence $\eps_n \searrow 0$ for which
the following holds: Let $K \subset \RR^n$ be a compact, convex
set with a non-empty interior. Let $X$ be a random vector that is
distributed uniformly in $K$. Then there exist a unit vector
$\theta$ in $\RR^n$, $t_0 \in \RR$ and $\sigma > 0$ such that
\begin{equation}
 \sup_{A \subset \RR} \left| \, Prob \, \{ \, \langle X, \theta
\rangle \in A \, \} \, - \, \frac{1}{\sqrt{2 \pi \sigma}} \int_A
e^{-\frac{(t - t_0)^2}{2 \sigma^2}} dt \, \right| \leq \eps_n,
\end{equation}
 where the supremum runs over all measurable sets $A
\subset \RR$, and where $\langle \cdot, \cdot \rangle$ denotes the
usual scalar product in $\RR^n$. Furthermore, under the additional
assumptions that the expectation of $X$ is zero and that the
covariance matrix of $X$ is the identity matrix, we may assert
that most unit vectors $\theta$ satisfy ($\ast$), with $t_0 = 0$
and $\sigma = 1$. Corresponding principles also hold for
multi-dimensional marginal distributions of convex sets.
\end{abstract}
\renewcommand{\theequation}{\arabic{equation}}

 \setcounter{equation}{0}

\section{Introduction}
\label{section1}
 We begin with an example. Let $n \geq 1$ be an
 integer, and consider the cube $Q = [ -\sqrt{3}, \sqrt{3}]^n
\subset \RR^n$. Suppose that $X = (X_1,...,X_n)$ is a random
vector that is distributed uniformly in the cube $Q$. Then
$X_1,...,X_n$ are independent, identically-distributed random
variables of mean zero and variance one. Consequently, the
classical central limit theorem states that the distribution of
the random variable
$$ \frac{X_1 + ... + X_n}{\sqrt{n}} $$
is close to the standard normal distribution, when $n$ is large.
Moreover, suppose we are given $\theta_1,...,\theta_n \in \RR$
with $\sum_{i=1}^n \theta_i^2 = 1$. Then
under mild conditions
on the $\theta_i$'s (such as Lindeberg's condition, see,
e.g., \cite[Section VIII.4]{feller}), the distribution of the random variable
$$ \langle \theta, X \rangle = \sum_{i=1}^n \theta_i X_i $$
is approximately gaussian, provided that
the dimension $n$ is large. For background on the classical
central limit theorem we refer the reader to, e.g., \cite{feller}
and \cite{stroock}.

\smallskip Let us consider a second example, no less fundamental
than the first. We denote by $| \cdot |$ the standard Euclidean
norm in $\RR^n$, and let $\sqrt{n + 2} \, D^{n} = \{ x \in \RR^n ;
|x| \leq \sqrt{n+2} \}$ be the Euclidean ball of radius
$\sqrt{n + 2}$ around the origin in $\RR^n$. We also write
$S^{n-1} = \{ x \in \RR^n ; |x| = 1 \}$ for the unit sphere in
$\RR^n$. Suppose that $Y = (Y_1,...,Y_n)$ is a random vector that
is distributed uniformly in the ball $\sqrt{n + 2} \, D^{n}$. Then
$Y_1,...,Y_n$ are identically-distributed random variables of mean
zero and variance one, yet they are not independent. Nevertheless,
it was already observed  by Maxwell that for any $\theta =
(\theta_1,...,\theta_n) \in
S^{n-1}$, the distribution of the random variable
$$ \langle \theta, Y \rangle = \sum_{i=1}^n \theta_i Y_i $$
is close to the standard normal distribution, when $n$ is large.
See, e.g., \cite{diaconis} for the history of the latter fact
and for more information.

\smallskip There is a wealth of central limit theorems in probability
theory that ensure normal approximation for a sum of many
independent, or weakly dependent, random variables. Our first
example, that of the cube, fits perfectly into this framework. The
approach we follow in this paper relates more to the second
example, that of the Euclidean ball, where the ``true source'' of
the gaussian approximation may be attributed to geometry. The
geometric condition we impose on the distribution of the random
variables  is that of convexity. We shall see that convexity may
substitute for independence in certain aspects of the phenomenon
represented by the classical central limit theorem.

\smallskip A function $f: \RR^n \rightarrow [0, \infty)$ is
log-concave if
$$ f(\lambda x + (1 - \lambda) y) \geq f(x)^{\lambda}
f(y)^{1-\lambda} $$ for all $x, y \in \RR^n$ and $0 < \lambda <
1$. That is, $f$ is log-concave when $\log f$ is concave on the
support of $f$. Examples of interest for log-concave functions
include characteristic functions of convex sets, the gaussian
density, and several densities from statistical mechanics.  In
this manuscript, we consider random vectors in $\RR^n$ that are
distributed according to a log-concave density. Thus, our
treatment includes as a special case the uniform distribution on
an arbitrary compact, convex set with a non-empty interior.

\smallskip We say that a function $f: \RR^n \rightarrow [0, \infty)$
is isotropic if it is the density of a random vector with zero
mean and identity covariance matrix. That is, $f$ is isotropic
when
$$ \int_{\RR^n} f(x) dx = 1, \ \ \ \int_{\RR^n} x f(x) dx = 0 \ \ \text{and} \ \ \
\int_{\RR^n} \langle x, \theta \rangle^2 f(x) dx = |\theta|^2
$$ for all $\theta \in \RR^n$. Any
log-concave function with $0 < \int f < \infty$ may be brought to
an isotropic position via an affine map, that is, $f \circ T$ is
isotropic for some affine map $T:\RR^n \rightarrow \RR^n$ (see,
e.g., \cite{MP}). Suppose that $X$ and $Y$ are two random
variables attaining values in some measure space $\Omega$ (here
$\Omega$ will always be $\RR$ or $\RR^n$ or a subspace $E \subset
\RR^n$). We define their total-variation distance as
$$ d_{TV}(X, Y) =
2  \sup_{A \subset \Omega} \left| \, Prob \left \{ X \in A \right
\} \, - \, Prob \left \{ Y \in A \right \} \, \right|,
$$ where the supremum runs over all measurable sets $A \subset
\Omega$. Note that $ d_{TV}(X, Y)$ equals the $L^1$-distance
between the densities of $X$ and $Y$, when these densities exist.
Let $\sigma_{n-1}$ stand for the unique rotationally-invariant
probability measure on $S^{n-1}$, also referred to as the uniform
probability measure on the sphere $S^{n-1}$.

\begin{theorem} There exist sequences $\eps_n \searrow 0,
\delta_n \searrow 0$ for which the following holds: Let $n \geq 1$,
and let $X$ be a random vector in $\RR^n$ with an isotropic, log-concave
density.
 Then there exists a subset $\Theta \subset S^{n-1}$
with $\sigma_{n-1}(\Theta) \geq 1 - \delta_n$, such that for all
$\theta \in \Theta$,
$$
d_{TV} \left( \, \langle X, \theta \rangle \, , \,
Z \, \right)
\leq \eps_n,
$$
where $Z \sim N(0,1)$ is a standard normal random variable.
\label{thm_basic}
\end{theorem}

We have the bounds $\eps_n \leq C \left( \frac{\log \log
(n+2)}{\log (n+1)} \right)^{1/2}$ and $\delta_n \leq \exp \left(
-c n^{0.99} \right)$ for $\eps_n$ and $\delta_n$ from Theorem
\ref{thm_basic}, where $c, C >0$ are universal constants. The
quantitative estimate we provide for $\eps_n$ is rather poor.
While Theorem \ref{thm_basic} seems to be a reasonable analog of
the classical central limit theorem for the category of
log-concave densities, we are still lacking the precise
Berry-Esseen type bound. A plausible guess might be that the
logarithmic dependence should be replaced by a power-type decay,
in the bound for $\eps_n$.

\smallskip Theorem \ref{thm_basic} implies the result stated in
the abstract of this paper, which does not require isotropicity;
indeed, recall that any log-concave density can be made isotropic
by applying an appropriate affine map. Thus, any log-concave
density in high dimension has at least one almost-gaussian
marginal. When the log-concave density is also isotropic, we can
assert that, in fact, the vast majority of its marginals are
approximately gaussian. An inherent feature of Theorem
\ref{thm_basic} is that it does not provide a specific unit vector
$\theta \in S^{n-1}$ for which $\langle X, \theta \rangle$ is
approximately normal.  This is inevitable: We clearly cannot take
$\theta = (1,0,...,0)$ in the example of the cube above, and hence
there is no fixed unit vector that suits all isotropic,
log-concave densities. Nevertheless, under additional symmetry
assumptions, we can identify a unit vector that always works.

\smallskip Borrowing terminology from Banach space theory, we say
that a function $f: \RR^n \rightarrow \RR$ is unconditional if
$$ f(x_1,...,x_n) = f(|x_1|,...,|x_n|) \ \ \ \text{for all} \ \ x =
(x_1,...,x_n) \in \RR^n. $$ That is, $f$ is unconditional when it
is invariant under coordinate reflections.

\begin{theorem} There exists a sequence $\eps_n \searrow 0$ for
which the following holds: Let $n \geq 1$, and let $f:\RR^n
\rightarrow [0, \infty)$ be an unconditional, isotropic,
log-concave function. Let $X = (X_1,...,X_n)$ be a random vector
in $\RR^n$ that is distributed according to the density $f$. Then,
$$ d_{TV} \left( \,
\frac{X_1 + ... + X_n}{\sqrt{n}} \, , \, Z \, \right) \leq \eps_n
$$ where $Z \sim N(0,1)$ is a standard normal random variable.
\label{thm_uncond}
\end{theorem}

We provide the estimate $\eps_n \leq \frac{C}{(\log (n+1))^{1/5}}$
for $\eps_n$ from Theorem \ref{thm_uncond}. Multi-dimensional
versions of Theorem \ref{thm_basic} are our next topic. For
integers $k,n$ with $1 \leq k \leq n$, let $G_{n,k}$ stand for the
grassmannian of all $k$-dimensional subspaces in $\RR^n$. Let
$\sigma_{n,k}$ be the unique rotationally-invariant probability
measure on $G_{n,k}$. Whenever we refer to the uniform measure on
$G_{n,k}$, and whenever we select a random $k$-dimensional
subspace in $\RR^n$, we always relate to the probability measure
$\sigma_{n,k}$ defined above. For a subspace $E \subset \RR^n$ and
a point $x \in \RR^n$, let
 $Proj_E(x)$ stand for the orthogonal projection of $x$ onto $E$.
A standard gaussian random vector in a $k$-dimensional subspace $E
\subset \RR^n$ is a random vector $X$ that satisfies $Prob \{ X
\in A \} = (2 \pi)^{-k/2} \int_A \exp(-|x|^2/2) dx$ for any
measurable set $A \subset E$.

\begin{theorem} There exists a universal constant $c > 0$
for which the following holds: Let $n \geq 3$ be an integer, and
let $X$ be a random vector in $\RR^n$ with an isotropic,
log-concave density. Let $\eps > 0$ and suppose that $1 \leq k
\leq c \eps^2 \frac{\log n}{\log \log n}$ is an integer. Then
there exists a subset $\E \subset G_{n,k}$ with $\sigma_{n,k}(\E)
\geq 1 - e^{-c n^{0.99}}$ such that for any $E \in \E$,
$$ d_{TV} \left( \, Proj_E(X) \, , \, Z_E \, \right) \leq \eps, $$
where $Z_E$ is a standard gaussian random vector in the subspace
$E$. \label{thm_multi}
\end{theorem}

That is, most $k$-dimensional marginals of an isotropic,
log-concave function, are approximately gaussian with respect to
the total-variation metric, provided that $k <<  \frac{\log
n}{\log \log n}$. Note the clear analogy between Theorem
\ref{thm_multi} and Milman's precise quantitative theory of
Dvoretzky's theorem, an analogy that dates back to Gromov
\cite[Section 1.2]{gromov}. Readers that are not familiar with
Dvoretzky's theorem are referred to, e.g., \cite[Section 4.2]{GM},
 to \cite{dvo_30_years} or to \cite{linden}. Dvoretzky's theorem shows that
$k$-dimensional geometric projections of an $n$-dimensional convex
body are $\eps$-close to a Euclidean ball, provided that $k <c
\eps^2 \log n$. Theorem \ref{thm_multi} states that
$k$-dimensional marginals, or measure-projections, of an
$n$-dimensional convex body are $\eps$-close to gaussian when $k <
c \eps^2 \log n / (\log \log n)$. Thus, according to Dvoretzky's
theorem, the geometric shape of the support of the marginal
distribution may be approximated by a very regular body -- a
Euclidean ball, or an ellipsoid -- whereas Theorem \ref{thm_multi}
demonstrates that the marginal distribution itself is very
regular;  it is approximately normal.

\smallskip More
parallels between Theorem \ref{thm_multi} and Dvoretzky's theorem
are apparent from the proof of Theorem \ref{thm_multi} below. We
currently do not know whether there exists a single subspace that
satisfies both the conclusion of Theorem \ref{thm_multi} and the
conclusion of Dvoretzky's theorem simultaneously; both theorems
show that a ``random subspace'' works with large probability, but
with respect to different Euclidean structures. The logarithmic
dependence on the dimension is known to be tight in Milman's form
of Dvoretzky's theorem. However, we have no reason to believe that
the quantitative estimates in Theorem \ref{thm_multi} are the best
possible.

\smallskip There are several mathematical articles
where Theorem \ref{thm_basic} is explicitly conjectured. Brehm and
Voigt suggest Theorem \ref{thm_basic} as a conjecture in
\cite{brehm_voigt}, where they write that this conjecture appears
to be ``known among specialists''. Anttila, Ball and Perissinaki
formulated the same conjecture in \cite{ABP}, independently and
almost simultaneously with Brehm and Voigt. Anttila, Ball and
Perissinaki also proved the conjecture for the case of uniform
distributions on  convex sets whose modulus of convexity and
diameter satisfy certain quantitative assumptions. Gromov wrote a
remark in \cite[Section 1.2]{gromov} that seems related  to
Theorem \ref{thm_basic} and Theorem \ref{thm_multi}, especially in
view of the techniques we use here. Following \cite{ABP} and
\cite{brehm_voigt}, significant contributions regarding the
central limit problem for convex sets were made by Bastero and
Bernu\'es \cite{bb}, Bobkov \cite{bobkov}, Bobkov and Koldobsky
\cite{bobkov_koldobsky}, Brehm and Voigt \cite{brehm_voigt},
 Brehm, Hinow, Vogt and Voigt \cite{BHVV},
Koldobsky and Lifshits \cite{koldobsky_lifshits}, E. and M. Meckes
\cite{meckes}, E. Milman \cite{emanuel}, Naor and Romik \cite{NR},
Paouris \cite{pa_clt}, Romik \cite{romik_phd}, S. Sodin
\cite{sodin}, Wojtaszczyk \cite{woj} and others.

\smallskip Let us explain a few ideas from our proof.
We begin with a general principle that goes back
to Sudakov \cite{sudakov} and to Diaconis and Freedman \cite{dia_thin_shell}
(see also the expositions of  Bobkov \cite{bobkov} and
 von Weizs\"acker \cite{weizs}.
A sharpening for the case of convex bodies was obtained by
Anttila, Ball and Perissinaki \cite{ABP}). This principle reads as
follows: Suppose $X$ is any random vector in $\RR^n$ with zero
mean and identity covariance matrix. Then most of the marginals of
$X$ are approximately gaussian, if and only if the random variable
$|X| / \sqrt{n}$ is concentrated around the value one. In other
words, typical marginals are approximately gaussian if and only if
most of the mass is concentrated on a ``thin spherical shell'' of
radius $\sqrt{n}$ and width much smaller than $\sqrt{n}$.
Therefore, to a certain extent, our task is essentially reduced to
proving the following:

\begin{theorem} Let $n \geq 1$ be an integer and let $X$ be
a random vector with an isotropic, log-concave density in $\RR^n$.
Then for all $0 \leq \eps \leq 1$,
$$
 Prob \left \{ \left| \, \frac{|X|}{\sqrt{n}} -1 \, \right| \geq \eps
\right \} \leq C n^{-c \eps^2}, $$ where $c, C > 0$ are universal
constants.
  \label{cor_202}
\end{theorem}

A significantly superior estimate to that of Theorem
\ref{cor_202}, for the case where $\eps$ is a certain universal
constant greater than one, is given by Paouris \cite{Pa2},
\cite{Pa3}. It would be interesting to try and improve the bound
in Theorem \ref{cor_202}  also for smaller values of $\eps$.

\smallskip Returning to the sketch of the proof,
suppose that we are given a random vector $X$ in $\RR^n$ with an
isotropic, log-concave density. We need to show that most of its
marginals are almost-gaussian. Select a random $k$-dimensional
subspace $E \subset \RR^n$, for a certain integer $k$. We use a
concentration of measure inequality -- in a way similar to
Milman's proof of Dvoretzky's theorem -- to show that with large
probability of choosing the subspace $E$, the distribution of the
random vector $Proj_E(X)$ is approximately spherically-symmetric.
This step is carried out in Section \ref{section3}, and it is also
outlined by Gromov \cite[Section 1.2]{gromov}.

\smallskip Fix a subspace $E$ such that $Proj_E(X)$
is approximately spherically-symmetric. In Section \ref{section4}
we use the Fourier transform to conclude that the approximation by
a spherically-symmetric distribution actually holds in the
stronger $L^{\infty}$-sense, after convolving with a gaussian. In
Section \ref{section5} we show that the gaussian convolution has
only a minor effect, and we obtain a spherically-symmetric
approximation to $Proj_E(X)$ in the total-variation, $L^1$-sense.
Thus, we obtain a density in the subspace $E$ that has two
properties: It is log-concave, by Pr\'ekopa-Leindler, and it is
also approximately radial. A key observation is that such
densities are necessarily very close to the uniform distribution
on the sphere; this observation boils down to estimating the
asymptotics of some one-dimensional integral. At this point, we
further project our density, that is already known to be close to
the uniform distribution on a sphere, to any lower-dimensional
subspace. By Maxwell's principle we obtain an approximately
gaussian distribution in this lower-dimensional subspace. This
completes the rough sketch of our proof.

\smallskip Throughout this paper, unless stated otherwise,
the letters
$c, C, c^{\prime}, \tilde{C}$ etc. denote positive
 universal constants, that are not necessarily the
same in different appearances. The symbols $C, C^{\prime},
\bar{C}, \tilde{C}$ etc. denote universal constants that are
assumed to be sufficiently large, while $c, c^{\prime}, \bar{c},
\tilde{c}$ etc. denote sufficiently small universal constants. We
abbreviate $\log$ for the natural logarithm, $\EE$ for
expectation, $Prob$ for probability and $Vol$ for volume.

\smallskip
\emph{Acknowledgements.} I would like to thank Charles Fefferman,
Emanuel Milman and Vitali Milman for interesting discussions on
related subjects, and to Boris Tsirelson for mentioning the
central limit problem for convex sets in his graduate course at
Tel-Aviv University.

 \setcounter{equation}{0}
\section{Some background on log-concave
functions} \label{section2}

 Here we gather some useful facts
pertaining mostly to log-concave densities. For more information
about log-concave functions, the reader is referred to, e.g.,
\cite{ball_log_concave}, \cite{dedicata} and \cite{lovasz}. The
raison d'\^{e}tre
 of log-concave densities on $\RR^n$ stems from the
classical Brunn-Minkowski inequality and its generalizations.
 Let $E \subset \RR^n$ be a
subspace, and let $f: \RR^n \rightarrow [0, \infty)$ be an
integrable function. We denote the marginal of $f$ with respect to
 the subspace $E$ by
$$ \pi_E(f)(x) = \int_{x + E^{\perp}} f(y) dy \ \ \ (x \in E) $$
where $x + E^{\perp}$ is the affine subspace in $\RR^n$ that is
orthogonal to $E$ and passes through $x$. The Pr\'{e}kopa-Leindler
inequality (see \cite{Pr1}, \cite{Le}, \cite{Pr2} or the first
pages of \cite{pisier_book}), which is a functional version of
Brunn-Minkowski, implies that  $\pi_E(f)$ is log-concave whenever
$f$ is log-concave and integrable. Therefore, when $f$ is
isotropic and log-concave, $\pi_E(f)$ is also isotropic and
log-concave. A further consequence of the Pr\'{e}kopa-Leindler
inequality, is that when $f$ and $g$ are integrable log-concave
functions on $\RR^n$, so is their convolution $f * g$. (The latter
result actually goes back to \cite{davidovic}, \cite{leker} and
\cite{schoenberg}.)

\begin{lemma} Let $n \geq 1$ be an integer,
and let $X$ be a random vector in $\RR^n$ with a log-concave
density. Assume that $F: \RR^n \rightarrow [0, \infty)$ is an
even, convex function, such that $F(t x) = t F(x)$ for all $t > 0,
x \in \RR^n$. Denote $E = \sqrt{\EE |F(X)|^2}$.
 Then,
\begin{enumerate}
\item[(i)] $\displaystyle
Prob \left \{ F(X) \geq t E  \right
\} \leq 2 e^{-t/10}$ for all $t \geq 0$.
\end{enumerate}
Additionally, let $0 < \eps \leq \frac{1}{2}$, and let $M > 0$
satisfy $ Prob \{ F(X) \geq M \} \leq \eps. $
Then,
\begin{enumerate}
\item[(ii)] $\displaystyle
Prob \left \{ F(X) \geq t M \right \} \leq (1 - \eps)
\left(\frac{\eps}{1-\eps} \right)^{(t + 1)/2}$ for all $t \geq 1$.
\end{enumerate}
\label{lem_1010}
\end{lemma}

Lemma \ref{lem_1010} is the well-known Borell's lemma (see its
elegant proof in \cite{borell_what} or \cite[Theorem III.3]{MS}).
Let $f: \RR^n \rightarrow [0, \infty)$ be an integrable function.
For $\theta \in S^{n-1}$ and $t \in \RR$ we define $H_{\theta, t}
= \{ x \in \RR^n ; \langle x, \theta \rangle \leq t \}$ and
\begin{equation}
M_f(\theta, t) = \int_{H_{\theta,t}} f(x) dx. \label{eq_146}
\end{equation}
The function $M_f$ is continuous in $\theta$ and $t$,
non-decreasing in $t$, and its derivative $\frac{\partial
M_f}{\partial t}$ is the Radon transform of $f$. Thus, in
principle, one may recover the function $f$ from a complete
knowledge of $M_f$. Clearly, for any subspace $E \subset \RR^n$,
\begin{equation}
 M_{\pi_E(f)}(\theta, t) = M_f(\theta, t) \ \ \ \text{for all} \ \ \ \theta \in S^{n-1} \cap E, t \in \RR.
\label{proj_ok}
\end{equation}
Moreover, let $\theta \in S^{n-1}$, let $E = \RR \theta$ be
the one-dimensional subspace spanned by $\theta$, and denote
 $g = \pi_E(f)$. Then
\begin{equation}
 g(t \theta) = \frac{\partial}{\partial t}
M_{\pi_E(f)} (\theta, t) = \frac{\partial}{\partial t} M_{f}
(\theta, t) \label{eq_621_}
\end{equation}
for all points $t \in \RR$ where, say, $g(t \theta)$ is
continuous.

\begin{lemma} Let $n \geq 1$ be an integer, and let
 $f: \RR^n \rightarrow [0, \infty)$ be an isotropic,
log-concave function. Fix $\theta \in S^{n-1}$. Then,
\begin{enumerate}
\item[(i)] For $t \geq 0$ we have  $
 1 - 2 e^{-|t|/10} \leq M_f(\theta, t) \leq 1$.
 \item[(ii)] For $t \leq 0$ we have $
0 \leq M_{f} \left( \theta, t \right) \leq 2 e^{-|t|/10}.$
\end{enumerate}
\label{lem_1028}
\end{lemma}

\emph{Proof:} Let $X$ be a random vector with density $f$. Then
$\EE |\langle X, \theta \rangle|^2 = 1$. We use Lemma
\ref{lem_1010}(i), with the function
 $F(x) = |\langle x,
\theta \rangle|$, to deduce the desired inequalities.
\hfill $\square$

\smallskip The space of all isotropic, log-concave functions in a
fixed dimension is a compact space, with respect to, e.g., the
$L^1$-metric. In particular, one-dimensional log-concave functions
are quite rigid. For instance, suppose that $g: \RR \rightarrow
[0, \infty)$ is an isotropic, log-concave function. Then (see
Hensley \cite{hensley} and also, e.g., \cite[Lemma 5.5]{lovasz} or
\cite{fradelizi}),
\begin{equation}
\frac{1}{10} \leq g(0) \leq \sup_{x \in \RR} g(x) \leq 1.
\label{eq_547}
\end{equation}
We conclude that for any log-concave, isotropic function $f: \RR^n
\rightarrow [0, \infty)$,
\begin{equation}
\left|  M_f(\theta, t)  - M_f(\theta, s)  \right| \leq |t - s| \ \
\ \text{for all} \ s, t \in \RR, \ \theta \in S^{n-1}.
\label{eq_547_}
\end{equation}
To prove (\ref{eq_547_}), we set $E = \RR \theta$ and $g =
\pi_E(f)$. Then $g$ is isotropic and log-concave, hence $\sup g
\leq 1$ by (\ref{eq_547}). Note that $g$ is continuous in the
interior of its support, since it is a log-concave function.
According to (\ref{eq_621_}), the function $t \mapsto g(t \theta)$
is the derivative of the function $t \mapsto M_{f} (\theta, t)$,
and (\ref{eq_547_}) follows.

\smallskip
Our next proposition is essentially taken from Anttila, Ball and
Perissinaki \cite{ABP}, yet we use the extension to the non-even
case which is a particular case of a result of Bobkov
\cite[Proposition 3.1]{bobkov}. A function $g: S^{n-1} \rightarrow
\RR$ is $L$-Lipshitz, for $L > 0$, if $|g(x) - g(y)| \leq L |x-y|$
for all $x, y \in S^{n-1}$.

\begin{proposition} Let $n \geq 1$ be an integer.
Let $t \in \RR$ and let $f: \RR^n \rightarrow [0, \infty)$ be an
isotropic, log-concave function. Then, the function
$$ \theta \mapsto M_f(\theta, t) \ \ \ (\theta \in S^{n-1}) $$
is $C$-Lipshitz on $S^{n-1}$. Here, $C > 0$ is a universal
constant. \label{lem_bobkov}
\end{proposition}

The proof of Proposition \ref{lem_bobkov} in \cite{bobkov}
involves analysis of two-dimensional log-concave functions. A
beautiful argument yielding Proposition \ref{lem_bobkov}, for the
case where $f$ is an even function, appears in \cite{ABP}. The
approach in \cite{ABP} is based on an application of Busemann's
theorem in dimension $n+1$, which leads to the conclusion  that
$\theta \mapsto |\theta| M_f(t, \theta / |\theta|)^{-1}$ is a norm
on $\RR^n$ for any fixed $t \geq 0$.

\setcounter{equation}{0}
\section{Techniques from Milman's proof of
Dvoretzky's theorem} \label{section3}

It is well-known that for large $n$, the uniform probability
measure $\sigma_{n-1}$ on the unit sphere $S^{n-1}$ satisfies
strong concentration inequalities. This concentration of measure
phenomenon is one of the main driving forces in high-dimensional
convex geometry, as was first demonstrated by Milman in his proof
of Dvoretzky's theorem
 (see \cite{mil70} or \cite[Section
4.2]{GM}).
 Our next proposition
is essentially taken from Milman's work, though the precise
formulation we use is due to Gordon \cite{gordon1}, \cite{gordon2}
(see also \cite{schechtman1}, \cite{schechtman2} or \cite[Theorem
6]{NR}).

\begin{proposition} Let $n \geq 1$ be an
integer, let $L > 0$, $0 < \eps \leq 1/2$, and let $g: S^{n-1}
\rightarrow \RR$ be an $L$-Lipshitz function. Denote $M =
\int_{S^{n-1}} g(x) d \sigma_{n-1}(x)$. Assume that $1 \leq k \leq
\hat{c} \eps^2 n$ is an integer. Suppose that
$E \in G_{n,k}$ is a random subspace, i.e., $E$ is distributed
according to the probability measure $\sigma_{n,k}$ on $G_{n,k}$.
Then, with probability greater than $1 - \exp \left(- c \eps^2 n
\right)$,
\begin{equation}
|g(\theta) - M| \leq \eps L \ \ \ \ \text{for all} \ \theta \in
S^{n-1} \cap E. \label{eq_700}
\end{equation} Here, $0 < c,\hat{c} < 1$ are universal constants.
\label{dvoretzky}
\end{proposition}

\smallskip Our use of ``Dvoretzky's theorem type'' arguments
 in the next lemma is inspired by the powerful
methods of Paouris in \cite{Pa1}, \cite{Pa2}, \cite{Pa3}.

\begin{lemma} Let $n \geq 1$ be an integer, let $A \geq 1, 0 < \delta \leq \frac{1}{2}$
and let $f: \RR^n \rightarrow [0, \infty)$ be an isotropic,
log-concave function. Assume that $1 \leq \ell \leq c \delta
A^{-1} \log n$ is an integer, and let $E$ be a random
$\ell$-dimensional subspace in $\RR^n$. Then with probability
greater than $1 - e^{-c n^{1-\delta}}$,
\begin{equation}
\sup_{\theta \in S^{n-1} \cap E} M_{f}(\theta, t) \leq e^{-A \ell}
+ \inf_{\theta \in S^{n-1} \cap E} M_{f}(\theta, t) \ \ \
\text{for all} \ \ t \in \RR. \label{eq_616}
\end{equation}
Here, $0 < c < 1$ is a universal constant.
\label{lem_dvo}
\end{lemma}

\emph{Proof:} We may assume that $n$ exceeds a given universal
constant, since otherwise, for a suitable choice of a small
universal constant $c$, there is no $\ell$ with $1 \leq
 \ell
\leq c \delta A^{-1} \log n$. Fix a real number $t$. According to
Proposition \ref{lem_bobkov}, the function $\theta \mapsto
M_f(\theta, t)$ is $C$-Lipshitz on $S^{n-1}$. Let $E \in G_{n,
\ell}$ be a random subspace, uniformly distributed in $G_{n,
\ell}$. We would like to apply Proposition \ref{dvoretzky} with $k
= \ell, L = C$ and $\eps = \frac{1}{2} n^{-\delta/2}$. Note that
for this choice of parameters,
$$ k = \ell \leq c \delta A^{-1} \log n \leq \hat{c} \eps^2 (\log
1/\eps)^2 n \ \ \ \text{and} \ \ \ 2 \eps L \leq e^{-2 c
\delta \log n} \leq e^{-2 A \ell},
$$ provided that $c$ is a sufficiently small, positive universal
constant, and that $n$ is greater than some universal constant.
Hence the appeal to Proposition \ref{dvoretzky} is legitimate.
From the conclusion of that proposition,  with probability larger
than $1 - e^{-c^{\prime} n^{1-\delta}}$ of selecting $E$,
\begin{equation}
\label{eq_434} \sup_{\theta \in S^{n-1} \cap E} M_f(\theta,t) \leq
 e^{-2 A \ell}  + \inf_{\theta \in S^{n-1} \cap E} M_f(\theta, t).
\end{equation}
For any fixed $t \in \RR$, the estimate (\ref{eq_434}) holds with
 probability greater than $1 - e^{-c^{\prime} n^{1-\delta}}$.
Denote $I = \{ i \cdot e^{-2 A \ell} \, ; \, i=-\lceil e^{30 A
\ell} \rceil,...,\lceil e^{30 A \ell} \rceil \}$.  Then, with
probability greater than $1 - e^{ -\bar{c} n^{1-\delta}}$, we
obtain
\begin{equation}
\forall t \in I, \ \ \ \sup_{\theta \in S^{n-1} \cap E} M_{f}
(\theta, t) \leq e^{-2 A \ell} +  \inf_{\theta \in S^{n-1} \cap E}
M_{f} (\theta, t). \label{eq_443}
\end{equation}
Indeed, the estimate for the probability follows from the
inequality $(2 e^{30 A \ell} + 3) e^{ - c^{\prime} n^{1-\delta}}
\leq e^{-\bar{c} n^{1-\delta}}$.

Fix an $\ell$-dimensional subspace $E
\subset \RR^n$ that
 satisfies
(\ref{eq_443}). Select $\theta_1, \theta_2 \in S^{n-1} \cap E$. We
will demonstrate that for any $t \in \RR$,
\begin{equation}
 M_{f} (\theta_1, t) \leq e^{-A \ell} +
M_{f} (\theta_2, t).
\label{eq_618}
\end{equation}
To that end, note that when $|t| \geq 20 A \ell$, by Lemma
\ref{lem_1028},
\begin{equation}
\left| M_f(\theta_1, t) - M_f (\theta_2, t) \right| \leq 2 e^{-|t| / 10}
\leq 2 e^{-2 A \ell} \leq e^{-A \ell}.
\label{eq_626}
\end{equation}
Hence (\ref{eq_618}) holds for $|t| \geq 20 A \ell$. We still need
to consider the case where $|t| < 20 A \ell$. In this case, $|t|
\leq e^{20 A \ell}$ and hence there exists $t_0 \in I$ with  $|t -
t_0| \leq \frac{1}{2} \cdot e^{-2 A \ell}$. According to
(\ref{eq_547_}) from Section \ref{section2}, the function $t
\mapsto M_{f} (\theta_i, t)$ is $1$-Lipshitz for $i=1,2$.
Therefore, by using (\ref{eq_443}), we conclude (\ref{eq_618})
also for the case where $|t| < 20 A \ell$. Thus (\ref{eq_618})
holds for all $t \in \RR$, under the assumption that $E$ satisfies
(\ref{eq_443}).

\smallskip
Recall that $\theta_1, \theta_2 \in S^{n-1} \cap E$ are arbitrary,
hence we may take the supremum over $\theta_1$ and the infimum
over $\theta_2$ in (\ref{eq_618}). We discover that whenever the
subspace $E$ satisfies (\ref{eq_443}), it necessarily also
satisfies (\ref{eq_616}). The probability for a random
$\ell$-dimensional subspace $E \subset \RR^n$ to satisfy
(\ref{eq_443}) was shown to be greater than $1 - e^{-\bar{c} n^{1
-\delta}}$. The lemma thus follows. \hfill $\square$

\smallskip
 \emph{Remark.}
For the case where $f$ is even, Lemma \ref{lem_dvo} follows
from a direct application of Dvoretzky's theorem in Milman's
form. Indeed, in this case, $\theta \mapsto |\theta| M_f(\theta,
t)^{-1}$ is a norm, and Lemma \ref{lem_dvo} asserts that this norm
is almost Hilbertian when restricted to certain random subspaces.

\setcounter{equation}{0}
\section{Almost spherical log-concave
functions} \label{section4}

A large portion of this section is devoted to proving the
following proposition.

\begin{proposition}
There exist universal constants $C_0,C > 1$ and $0 < c < 1$ for
which the following holds: Let $n \geq 1$ be an integer and let
$f: \RR^n \rightarrow [0, \infty)$ be an isotropic, log-concave
function. Assume that
\begin{equation}
\sup_{\theta \in S^{n-1}} M_{f}(\theta, t) \leq e^{-C_0 n} +
\inf_{\theta \in S^{n-1}} M_{f}(\theta, t) \ \ \ \text{for all} \
\ t \in \RR. \label{eq_848_}
\end{equation}
Suppose that $Y$ is a random vector in $\RR^n$
with density $f$.
Then for all $0 < \eps < 1$,
\begin{equation}
 Prob \left \{ \left| \, \frac{|Y|}{\sqrt{n}} -1 \, \right| \geq \eps
\right \} \leq C e^{-c \eps^2 n }. \label{eq_309}
\end{equation}   \label{cor_735}
\end{proposition}

For $n \geq 1$ and $v > 0$ we define $\gamma_{n,v}: \RR^n
\rightarrow [0, \infty)$ to be the function \begin{equation}
\gamma_{n,v}(x) = \frac{1}{(2 \pi v)^{n/2}} \exp \left(
-\frac{|x|^2}{2v} \right). \label{gamma_v}
\end{equation}
 Then
$\gamma_{n,v}$ is the density of a gaussian random vector in
$\RR^n$ with expectation zero and covariance  matrix that equals
$v Id$, where $Id$ is the identity matrix. We write $O(n)$ for the
group of orthogonal transformations of $\RR^n$.

\begin{lemma} Let $n \geq 1$ be an integer,
let $\alpha \geq 5$, and let $f: \RR^n \rightarrow [0, \infty)$ be
an isotropic, log-concave function. Assume that
\begin{equation}
\sup_{\theta \in S^{n-1}} M_{f}(\theta, t) \leq e^{-5 \alpha n} +
\inf_{\theta \in S^{n-1}} M_{f}(\theta, t) \ \ \ \text{for all} \
\ t \in \RR. \label{eq_848}
\end{equation}
Denote $g = f * \gamma_{n, 1}$, where $*$ stands for convolution.
Then,
\begin{equation}
\sup_{\theta \in S^{n-1}} g(t \theta) \leq e^{-\alpha n} +
\inf_{\theta \in S^{n-1}} g(t \theta) \ \ \ \text{for all} \ \  t
\geq 0. \label{eq_750} \end{equation} \label{lem_fourier}
\end{lemma}

\emph{Proof:} We will show that the Fourier transform of $f$ is
almost spherically-symmetric. As usual, we define
$$
 \hat{f}(\xi) = \int_{\RR^n} e^{-2 \pi i \langle \xi, x \rangle}
f(x) dx \ \ \ \ (\xi \in \RR^n), $$ where $i^2 = -1$. Let $r > 0$,
and fix $\xi_1, \xi_2 \in \RR^n$ with $|\xi_1| = |\xi_2| = r$.
Denote by $E_1 = \RR \xi_1, E_2 = \RR \xi_2$ the one-dimensional
subspaces spanned by $\xi_1, \xi_2$, respectively. From
(\ref{eq_621_}) of Section \ref{section2} we see that
$\pi_{E_j}(f)(t \xi_j/|\xi_j|) = \frac{\partial}{\partial t} M_{f}
(\xi_j / |\xi_j|, t)$ for $j=1,2$ and for all $t$ in the interior
of the support of the log-concave function $t \mapsto
\pi_{E_j}(f)(t \xi_j/|\xi_j|)$. By integrating by parts we obtain
\begin{eqnarray}
\nonumber \label{f_xi} \lefteqn{\hat{f}(\xi_1) - \hat{f}(\xi_2)=
\int_{-\infty}^{\infty} \left[ \pi_{E_1}(f) \left(t
\frac{\xi_1}{|\xi_1|} \right) - \pi_{E_2}(f) \left(t
\frac{\xi_2}{|\xi_2|} \right) \right] e^{-2 \pi i r t} dt } \\ & =
& 2 \pi i r \int_{-\infty}^{\infty} \left[M_{f} \left(
\frac{\xi_1}{|\xi_1|}, t \right) - M_{f} \left(
\frac{\xi_2}{|\xi_2|}, t \right)  \right] e^{-2 \pi i r t} dt,
\phantom{aaaaaaaa}
\end{eqnarray}
as the boundary terms clearly vanish. From Lemma \ref{lem_1028} we
have
\begin{equation}
\left| M_{f} \left( \frac{\xi_1}{|\xi_1|}, t \right) - M_{f}
\left( \frac{\xi_2}{|\xi_2|}, t \right)  \right| \leq 2
e^{-|t|/10} \ \ \ \text{for all} \ \ t \in \RR. \label{eq_933}
\end{equation}
According to (\ref{f_xi}), (\ref{eq_933}) and to our assumption
(\ref{eq_848}), we conclude that for any $r > 0$ and $\xi_1, \xi_2
\in \RR^n$ with $|\xi_1| = |\xi_2| = r$,
\begin{eqnarray}
 \label{eq_1042}
\lefteqn{|\hat{f}(\xi_1) - \hat{f}(\xi_2)| } \\ & \leq & 2 \pi r
\left[ 80 \alpha n  \cdot e^{-5 \alpha n} + \int_{|t|
> 40 \alpha n} 2 e^{-|t|/10} dt \right] \leq r e^{-2 \alpha n},
\nonumber
\end{eqnarray}
where we made use of the fact that $\alpha n \geq 5$. A standard
computation  (e.g. \cite[page 6]{stein_weiss}) shows that
$\widehat{\gamma_{n,1}}(\xi) = e^{-2 \pi^2 |\xi|^2}$. Recall that
we define $g = f * \gamma_{n,1}$, and hence $\hat{g} (\xi) = e^{-2
\pi^2 |\xi|^2} \cdot \hat{f} (\xi)$. We thus deduce from
(\ref{eq_1042}) that for any $\xi_1, \xi_2 \in \RR^n$,
\begin{equation}
|\hat{g}(\xi_1) - \hat{g}(\xi_2)| \leq   e^{-2 \pi^2 r^2} r e^{-2
\alpha n} \ \ \ \text{whenever} \ \ |\xi_1| = |\xi_2| = r > 0.
\label{eq_1103}
\end{equation}
Let $x \in \RR^n$, and let $U \in O(n)$ be an orthogonal
transformation. By using the inverse Fourier transform (see, e.g.
\cite[Chapter I]{stein_weiss}) and applying (\ref{eq_1103}), we
get
\begin{eqnarray}
\label{eqn_1110} \lefteqn{  | g(x) - g(U x) | =
 \nonumber  \left| \int_{\RR^n} \left[ \hat{g}(\xi) -
\hat{g}(U \xi) \right] e^{2 \pi i \langle x, \xi \rangle} d\xi
\right|} \\ & \leq &  \int_{\RR^n} e^{-2 \pi^2 |\xi|^2} |\xi|
e^{-2 \alpha n}  d \xi \leq e^{-2 \alpha n} \int_{\RR^n} e^{-\pi
|\xi|^2} d \xi = e^{-2 \alpha n}.
\end{eqnarray}
Since $x \in \RR^n$ and $U \in O(n)$ are arbitrary, from
(\ref{eqn_1110}) we conclude (\ref{eq_750}). \hfill $\square$

\smallskip Let $f: [0, \infty) \rightarrow
[0, \infty)$ be a log-concave function with $0 < \int_0^{\infty} f
< \infty$, that is continuous on $[0, \infty)$ and $C^2$-smooth on
$(0, \infty)$.
 For $p
> 1$, denote by $t_p(f)$ the unique $t > 0$ for which $f(t) > 0$
and also
\begin{equation}
(\log  f)^{\prime}(t) = \frac{f^{\prime}(t)}{f(t)} =
-\frac{p-1}{t}. \label{eq_106}
\end{equation}

\begin{lemma} $t_p(f)$ is well-defined, under the above
assumptions on $f$ and $p$.
\end{lemma}

\emph{Proof:} We need to explain why a solution $t$ to equation
(\ref{eq_106})
 exists and is unique, for all $p > 1$. To that end,
note that $f$ is a log-concave function with finite, positive
mass, hence it decays  exponentially fast at infinity (this is a
very simple fact; see, e.g., \cite[Lemma 2.1]{psitwo}). Therefore,
the function $\vphi(t) = t^{p-1} f(t)$ satisfies
$$ \lim_{t \rightarrow 0^+} \vphi(t) = \lim_{t \rightarrow \infty}
\vphi(t) = 0. $$ The function $\vphi$ is continuous, non-negative,
not identically zero, and tends to zero at $0$ and at $\infty$.
Consequently, $\vphi$ attains its positive maximum at some finite
point $t_0
> 0$. Then $\vphi(t_0) > 0$ and $\vphi^{\prime}(t_0) = 0$, since $\vphi$ is $C^2$-smooth.
  On the other hand, $f$ is log-concave, and $t \mapsto t^{p-1}$
is strictly log-concave, hence $\vphi$ is strictly log-concave on
its support. Therefore, there is at most one point where $\vphi$
is non-zero and $\vphi^{\prime}$ vanishes. We conclude that there
exists exactly one point $t_0 > 0$ such that $f(t_0) > 0$ and
$$ \vphi^{\prime}(t_0) = t_0^{p-2} \left[ (p-1) f(t_0) + t_0
f^{\prime}(t_0) \right] = 0. $$ Thus a finite, positive $t$ that
solves (\ref{eq_106}) exists and is unique. \hfill $\square$

\smallskip Let us mention a few immediate properties of the quantity $t_p(f)$.
First, $f(t_p(f)) > 0$ for all $p > 1$. Second, suppose that $f$
is a continuous, log-concave function on $[0, \infty)$,
$C^2$-smooth on $(0, \infty)$, with $0 < \int f < \infty$. Then,
\begin{equation}
f( t ) \geq e^{-(n-1)} f(0) \ \ \ \text{for any }  \ \  0 \leq t
\leq t_n(f). \label{eq_1127}
\end{equation}
Indeed, if $f(0) = 0$ then (\ref{eq_1127}) is trivial. Otherwise,
$f(0) > 0$ and $f(t_n(f)) > 0$, hence $f$ is necessarily positive
on $[0, t_n(f)]$ by log-concavity. Therefore $\log f$ is finite
and continuous on $[0, t_n(f)]$, and $C^2$-smooth in $(0,
t_n(f))$. Additionally, $\log f$ is concave,
 hence $(\log f)^{\prime}$ is non-increasing in $(0, t_n(f))$.
From the definition (\ref{eq_106}) we deduce that $(\log
f)^{\prime}(t) \geq -(n-1) / t_n(f)$ for all $0 < t < t_n(f)$, and
(\ref{eq_1127}) follows.

\smallskip
Furthermore, since $(\log f)^{\prime}$ is non-increasing on the
interval in which it is defined, then $(\log f)^{\prime}(t) \leq
-(n-1) / t_n(f)$ for $t > t_n(f)$ for which $f(t) > 0$. We
conclude that for any $\alpha \geq 1$,
\begin{equation}
f( t ) \leq e^{-(\alpha - 1) (n-1)} f(t_n(f)) \ \ \ \text{when}  \
\ t \geq  \alpha t_n(f). \label{eq_1127_2}
\end{equation}

\smallskip Note that $t_p(f)$
behaves well under scaling of $f$. Indeed, let $f$ be a
continuous, log-concave function on $[0, \infty)$, $C^2$-smooth on
$(0, \infty)$, with $0 < \int f < \infty$. For $\delta
> 0$, denote $\tau_{\delta}(x) = \delta x$. From the definition
(\ref{eq_106}) we see that for any $p > 1$,
\begin{equation}
t_p( f \circ \tau_{\delta} ) = \delta^{-1} \cdot t_p(f).
\label{scale}
\end{equation}

\begin{lemma} Let $n \geq 2$, and let $f,g : [0, \infty) \rightarrow [0,
\infty)$ be continuous, log-concave functions, $C^2$-smooth on
$(0, \infty)$, such that $f(0) > 0, g(0) > 0$ and $\int f <
\infty, \int g < \infty $. Assume that for any $t \geq 0$,
\begin{equation}
 |f(t) -  g(t)| \leq e^{-5 n} \min \{ f(0), g(0) \}.
\label{eq_1132_}
\end{equation}
Then,
$$ \left(1 - e^{-n} \right) t_n(g) \leq t_n(f) \leq \left(1 + e^{-n} \right)
t_n(g). $$ \label{lem_257}
\end{lemma}

\emph{Proof:} Set $\delta = t_n(f)$. According to (\ref{scale}),
both the conclusions and the requirements of the lemma are
invariant when we replace $f$ and $g$ with $f \circ \tau_{\delta}$
and $g \circ \tau_{\delta}$, respectively. We apply this
replacement, and assume from now on that $t_n(f) = 1$.

\smallskip Inequality
(\ref{eq_1127}) and our assumption that $f(0) > 0$ show that $f(t)
\geq e^{-n} f(0)
> 0$ for $0 \leq t \leq 1$. We combine this inequality with
(\ref{eq_1132_}) to obtain the bound $|g(t) / f(t) - 1| \leq e^{-4
n}$ for all $0 \leq t \leq 1$. In particular,
 $g$ is positive on
$[0, 1]$. Denote $f_0 = \log f, g_0 = \log g$. Then
for all $0 \leq t \leq 1$,
\begin{equation}
-2 e^{-4 n} < \log(1 - e^{-4 n}) \leq g_0(t) - f_0(t) \leq \log (1
+ e^{-4 n}) < e^{-4 n}. \label{eq_1135}
\end{equation}
Next, we claim that
\begin{equation}
g_0^{\prime}(t) \geq f_0^{\prime} \left(t + e^{-2 n} \right) - 4
e^{-2 n} \ \ \ \text{for all} \ \ 0 < t \leq 1 - e^{-2 n}.
\label{eq_1129}
\end{equation}
Indeed, assume by contradiction that (\ref{eq_1129}) does not
hold. Then there exists $0 < t_0 \leq 1 - e^{-2 n}$ for which
$g_0^{\prime}(t_0) < f_0^{\prime}(t_0 + e^{-2n}) - 4 e^{-2n}$.
From our assumptions, $f$ and $g$ are log-concave, hence $f_0$ and
$g_0$ are concave, and hence $f_0^{\prime}$ and $g_0^{\prime}$ are
non-increasing on $(0,1)$. Therefore, for $t \in (t_0, t_0 +
e^{-2n})$,
\begin{equation}
g_0^{\prime}(t) \leq g_0^{\prime}(t_0) < f_0^{\prime}(t_0 +
e^{-2n}) - 4 e^{-2n} \leq f_0^{\prime}(t) - 4 e^{-2n}.
\label{eq_1149}
\end{equation}
Denote $t_1 = t_0 + e^{-2n}$. Then $[t_0, t_1] \subset [0, 1]$ and
by (\ref{eq_1149}),
$$ \left[ f_0(t_1) - g_0(t_1) \right] \, - \, \left [f_0(t_0) - g_0(t_0) \right]
\, > \, 4 e^{-2n} \cdot (t_1 - t_0) \, = \, 4 e^{-4n},
$$ in contradiction to (\ref{eq_1135}). Thus, our momentary
assumption -- that (\ref{eq_1129}) does not hold -- was false, and
hence (\ref{eq_1129}) is proved.

\smallskip From the definition (\ref{eq_106})
we see that $f_0^{\prime}(1) = (\log f)^{\prime}(1) = -(n-1)$.
Recall once again that $g_0^{\prime}$ is non-increasing. By
applying the case $t = 1 - e^{-2n}$ in (\ref{eq_1129}), we
conclude that for  $0 < s < 1 - 4 e^{-2n}$,
\begin{eqnarray}
 \nonumber
\lefteqn{ g_0^{\prime}(s) \geq g_0^{\prime}(1 -
e^{-2n}) \geq
f_0^{\prime}(1) - 4 e^{-2n} =  -(n-1) - 4 e^{-2n} } \\
& \geq &  -(n-1) \left( 1 + 4 e^{-2n} \right)
\geq -\frac{n-1}{1 - 4 e^{-2n}} > -\frac{n-1}{s}. \phantom{aaaaaaaaa}
\label{eq_1135_}
\end{eqnarray}
From (\ref{eq_1135_}) we conclude that $g^{\prime}(s) / g(s) =
g_0^{\prime}(s) \neq -\frac{n-1}{s}$ for all $0 < s < 1 - 4
e^{-2n}$. The definition (\ref{eq_106}) shows that $$ t_n(g) \geq
1 - 4 e^{-2n}.
$$ Recalling the scaling argument above, we see that we have actually
proved that
$$ t_n(g) \geq (1 - 4 e^{-2n}) t_n(f), $$
whenever the assumptions of the lemma hold. However, these
assumptions are symmetric in $f$ and $g$. Hence,
$$ t_n(g) \geq (1 - 4 e^{-2n}) t_n(f) \ \ \ \text{and also} \ \ \ t_n(f) \geq (1 - 4 e^{-2n}) t_n(g)
$$
for any functions $f,g$ that satisfy the assumptions of the lemma.
Since $1 + e^{-n} \geq 1/(1 - 4 e^{-2n})$ for $n \geq 2$, the
lemma is proved. \hfill $\square$

\smallskip Our next lemma is a standard application
of the Laplace asymptotic method, and is similar to, e.g.,
 \cite[Lemma 2.1]{dedicata} and
\cite[Lemma 2.5]{psitwo}. We will make use of the following
well-known bound: For $\alpha, \delta > 0$,
\begin{equation}
\label{eq_1218} \int_{\delta}^{\infty} e^{-\alpha \frac{t^2}{2}}
dt = \frac{1}{\sqrt{\alpha}} \int_{\delta \sqrt{\alpha}}^{\infty}
e^{-\frac{t^2}{2}} dt \leq \frac{\sqrt{2 \pi}}{\sqrt{\alpha}}
e^{-\alpha \frac{\delta^2}{2}}.
\end{equation}
The inequality in (\ref{eq_1218}) may be proved, for example, by
computing the Laplace transform of the gaussian density and
applying Markov's inequality (e.g., \cite[Section 1.3]{stroock}).

\begin{lemma} Let $n \geq 2$ be an integer, and let $f: [0, \infty) \rightarrow [0, \infty)$
be a continuous, log-concave function, $C^2$-smooth on $(0,
\infty)$, with $0 < \int_0^{\infty} f < \infty$. Then for $0 \leq
\eps \leq 1$,
\begin{equation}
\int_{t_n(f) (1 - \eps)}^{t_n(f) (1 + \eps)} t^{n-1} f(t) dt \geq
\left(1 - C e^{-c \eps^2 n} \right) \int_{0}^{\infty} t^{n-1} f(t)
dt, \label{eq_248}
\end{equation}
 where $C > 1$ and $0 < c < 1$
are universal constants. \label{lem_1153}
\end{lemma}

\emph{Proof:} We begin with a scaling argument.  A glance at
(\ref{scale}) and (\ref{eq_248}) assures us that both the validity
of the assumptions and the validity of the conclusions of the
present lemma, are not altered when we replace $f$ with $f \circ
\tau_{\delta}$, for any $\delta > 0$. Hence, we may switch from
$f$ to $f \circ \tau_{t_n(f)}$, and reduce matters to the case
$t_n(f) = 1$. Thus $f(1) > 0$. Multiplying $f$ by an appropriate
positive constant, we may assume that $f(1) = 1$.

\smallskip
We denote  $\psi(t) = (n-1) \log t + \log f(t) \ \ (t > 0)$, where we set
 $\psi(t) = -\infty$ whenever $f(t) = 0$. Since $f(1) = 1$, then
 $\psi(1) = 0$. Additionally, $\psi^{\prime}(1) = 0$ because
 $t_n(f) = 1$.
The function $\psi$ is concave, and therefore it attains its
maximum at $1$. Let $s_0, s_1 > 0$ be the minimal positive numbers
for which $\psi(1 - s_0) = -1$ and $\psi(1 + s_1) = -1$. Such
$s_0$ and $s_1$ exist since $\psi$ is continuous, $\psi(1) = 0$
and $\psi(t) \rightarrow -\infty$ when $t \rightarrow 0$ (because
of $\log t$) and when $t \rightarrow \infty$ (because of $\log f$,
since $f$ is log-concave with $0 < \int f < \infty$).

\smallskip We may suppose that $n \geq 100$; for an appropriate choice of a
large universal constant $C$, the right hand side of
(\ref{eq_248}) is negative for $n < 100$, and hence the lemma is
obvious for $n < 100$. Denote $m = \inf \{ t > 0; f(t) \neq 0 \}$
and $M = \sup \{ t > 0; f(t) \neq 0 \}$. Since $t_n(f) = 1$,
necessarily $m < 1$ and $M > 1$. Then, for $m < t < M$,
\begin{equation}
\psi^{\prime \prime}(t) = -\frac{n-1}{t^2} + (\log f)^{\prime
\prime}(t) \leq -\frac{n-1}{t^2}, \label{eq_113}
\end{equation}
since $\log f$ is concave and hence $(\log f)^{\prime \prime} \leq
0$. From (\ref{eq_113}) we obtain, in particular, the inequality
$\psi^{\prime \prime}(t) \leq -\frac{n-1}{4}$ for $m < t < \min \{
2, M \}$. Recalling that $\psi(1) = \psi^{\prime}(1) = 0$, we see
that $\psi(t) \leq -\frac{n-1}{8}  (t -1)^2$ for all $0 < t < 2$.
Therefore $\psi(1 - 4 / \sqrt{n}) \leq -1$ and $ \psi(1 + 4 /
\sqrt{n}) \leq -1$, and consequently
\begin{equation} s_0 \leq \frac{4}{\sqrt{n}} \ \ \ \text{and} \ \ \ s_1 \leq \frac{4}{\sqrt{n}}.
\label{s0_above}
\end{equation}
Since $n \geq 100$, then (\ref{s0_above}) implies that $s_0, s_1
\leq \frac{1}{2}$. Recall that the function $\psi$ is concave,
hence $\psi^{\prime}$ is non-increasing. The relations $\psi(1 -
s_0) = \psi(1 + s_1) = -1, \psi(1) = 0$ thus imply that
\begin{equation}
\psi^{\prime}(1 - s_0) \geq \frac{1}{s_0} \ \ \ \text{and} \ \ \
\psi^{\prime}(1 + s_1) \leq -\frac{1}{s_1}. \label{eq_302}
\end{equation}
Examination of (\ref{eq_113}) shows us that $\psi^{\prime
\prime}(t) \leq -(n-1)$ for $m < t \leq 1 - s_0$. By definition,
$\psi(1 - s_0) = -1$. We thus conclude from (\ref{eq_302}) that
$\psi(1 - s_0 - t) \leq -1 -\frac{t}{s_0}- \frac{n-1}{2} t^2 \ $
for $0 < t < 1 - s_0$. Fix $0 \leq \eps \leq 1$. Then,
\begin{eqnarray}
\label{eq_121} \lefteqn{ \int_0^{1 - s_0 - \eps} e^{\psi(t)} dt
\leq e^{-1}  \int_{\eps}^{\infty} e^{-\frac{t}{s_0} -(n-1) \frac{t^2}{2}} dt
} \\ & \leq &
\nonumber
 \min \left \{ s_0 e^{-\frac{\eps}{s_0}}, \int_{\eps}^{\infty}
e^{-(n-1) \frac{t^2}{2}} dt \right \} \leq \min \left \{ s_0
e^{-\frac{\eps}{s_0}}, \frac{e^{-(n-1)
\frac{\eps^2}{2}}}{\sqrt{(n-1) / (2 \pi)}} \right \}
\end{eqnarray}
where we used (\ref{eq_1218}) to estimate  the last integral.
Next, observe again that  $\psi^{\prime \prime}(t) \leq -
\frac{n-1}{4}$ for all $m < t < \min \{ 2, M \}$, by
(\ref{eq_113}). We use (\ref{eq_302}), as well as the fact that
$\psi(1 + s_1) = -1$, to obtain
\begin{equation} \psi(1 + s_1 + t) \leq -1 -\frac{t}{s_1} - \frac{n-1}{8}
t^2 \ \ \ \ \ \text{for} \ \ 0 \leq t \leq 1 - s_1.
\label{eq_1046}
\end{equation}
Consequently,
\begin{eqnarray}
\label{eq_909}
\lefteqn{ \int_{1 + s_1 + \eps}^{2} e^{\psi(t)} dt
\leq e^{-1} \int_{\eps}^{\infty} e^{-\frac{t}{s_1} -(n-1) \frac{t^2}{8}} dt
} \\
 & \leq & \nonumber
 \min \left \{ s_1 e^{-\frac{\eps}{s_1}}, \int_{\eps}^{\infty}
e^{-(n-1) \frac{t^2}{8}} dt \right \} \leq \min \left \{ s_1
e^{-\frac{\eps}{s_1}}, \frac{e^{-(n-1) \frac{\eps^2}{8}}}{ \sqrt{
(n-1) / (8 \pi)}} \right \}
\end{eqnarray}
by (\ref{eq_1218}). Since $s_1 \leq \frac{1}{2}$, we
deduce from (\ref{eq_1046})  that $\psi(2) \leq  - \frac{1}{2s_1}
-\frac{n-1}{32}$. Recall that $\psi^{\prime}$ is non-increasing,
that $\psi^{\prime}(1) = 0$ and that $\psi^{\prime \prime}(t) \leq
- \frac{n-1}{4}$ for $1 < t < \min \{2, M \}$. Therefore,
$\psi^{\prime}(t) \leq -\frac{n-1}{4}$ whenever $2 \leq  t < M$.
Thus we realize that
 $\psi(2+t) \leq \left( -\frac{1}{2s_1} - \frac{n-1}{32} \right) - \frac{n-1}{4} t$ for
$t \geq 0$. Hence,
\begin{equation}
\int_{2}^{\infty} e^{\psi(t)} dt \leq e^{-\frac{1}{2s_1} -
\frac{n-1}{32}}
 \int_{0}^{\infty} e^{-\frac{n-1}{4} t} dt \leq \frac{8 s_1}{n-1} e^{-\frac{n-1}{32}}. \label{eq_935}
\end{equation}
Let $s = s_0 + s_1$. Then, by the definition of $s_0$ and $s_1$,
\begin{equation}
\int_0^{\infty} e^{\psi(t)} dt \geq \int_{1 - s_0}^{1 + s_1}
e^{\psi(t)} dt \geq  \int_{1 - s_0}^{1 + s_1} e^{-1} dt = e^{-1}
s. \label{total}
\end{equation}
The inequalities we gathered above will allow us to prove
(\ref{eq_248}). Note that (\ref{eq_248}) is trivial when $\eps
\leq \frac{4}{\sqrt{n}}$; for an appropriate choice of a large
constant $C$, the right-hand side of (\ref{eq_248}) is negative in
this case. We may thus restrict our attention to the case where
$\frac{4}{\sqrt{n}} < \eps < 1$. Hence, $s_0 + \eps \leq 2 \eps$
and $s_1 + \eps \leq 2 \eps$, by (\ref{s0_above}). We add
(\ref{eq_121}), (\ref{eq_909}) and (\ref{eq_935}) to get
\begin{equation}
\int_{|t - 1| \geq 2 \eps} e^{\psi(t)} dt \leq \min \left \{ s e^{-\eps / s},
 \frac{20}{\sqrt{n}}
e^{-\frac{\eps^2 n}{20}} \right \} + \frac{20 s }{n} \cdot
e^{-n/100}. \label{eqn_603}
\end{equation}
Division of (\ref{eqn_603}) by (\ref{total}) yields,
\begin{equation}
\frac{\int_{|t - 1| \geq 2 \eps} \exp (\psi(t))
dt}{\int_0^{\infty} \exp (\psi(t)) dt} \leq 60 \min \left \{
e^{-\eps/s}, \frac{e^{-\frac{\eps^2 n}{20}}}{s
\sqrt{n}}
 \right \} + 40 e^{-n/100}.
 \label{eq_613}
\end{equation}
In order to establish (\ref{eq_248}) and complete the proof,
it is sufficient to show that
\begin{equation}
\int_{|t - 1| \geq 2 \eps} \exp (\psi(t)) dt \leq  100 e^{-\eps^2
n / 100}  \int_0^{\infty} \exp (\psi(t)) dt. \label{eq_624}
\end{equation}
According to (\ref{s0_above}), we know
that $s = s_0 + s_1 \leq \frac{10}{\sqrt{n}}$. In the case where
$$ \eps > 10 \frac{\sqrt{\log \frac{10}{s
\sqrt{n}}}}{\sqrt{n}}, $$ we have $\frac{1}{s \sqrt{n}}
< \exp \left( \frac{\eps^2 n}{100} \right)$ and hence the estimate
(\ref{eq_624}) follows from (\ref{eq_613}) by choosing the
``$\frac{e^{-\frac{\eps^2 n}{20}}}{s \sqrt{n}}$'' term
in the minimum in (\ref{eq_613}). In the complementary case, we
have $$ \eps \leq 10 \frac{\sqrt{\log \frac{10}{s
\sqrt{n}}}}{\sqrt{n}} \leq \frac{100}{s n}, $$ since
$\sqrt{\log t} \leq t$ for $t \geq 1$.
In this case, $\eps / s \geq \frac{1}{100} \eps^2 n$, and (\ref{eq_624}) follows by
selecting the ``$e^{-\eps/s}$'' term in (\ref{eq_613}).
Hence (\ref{eq_624}) is proved for all cases. The proof is complete.
 \hfill $\square$

\smallskip The following lemma is standard, and is almost
identical, for example, to \cite[Appendix V.4]{MS}. For a random
vector $X$ in $\RR^n$, we denote its covariance matrix by
$Cov(X)$.

\begin{lemma} Let $n \geq 1$ be an integer, let $A, r, \alpha, \beta > 0$ and
let $X$ be a random vector in $\RR^n$ with $\EE X = 0$ and $Cov(X)
= \beta Id$. Assume that the density of $X$ is log-concave, and
that
\begin{equation}
 Prob \left \{ \left| \, \frac{|X|}{r} - 1 \, \right|  \geq \eps
\right \} \leq A e^{-\alpha \eps^2 n } \ \ \ \ \text{for} \  \ 0
\leq \eps \leq 1. \label{eq_855}
\end{equation}
Then,
\begin{enumerate}
\item[(i)] For all $\displaystyle
\  0 \leq \eps \leq 1, \ \ \ \
 Prob \left \{  \left| \, \frac{|X|}{\sqrt{\beta n}} - 1 \, \right|  \geq \eps
\right \} \leq C^{\prime} e^{-c^{\prime} \eps^2 n }. $
\item[(ii)] $\displaystyle  \left| \frac{r}{\sqrt{\beta n}} - 1
\right| \leq \frac{C}{\sqrt{n}} \ $ provided that $n \geq C$.
\end{enumerate}
Here, $C, C^{\prime}, c^{\prime} > 0$
are constants that depend solely on $A$ and $\alpha$. \label{crude_lem}
\end{lemma}

\emph{Proof:} By a simple scaling argument, we may assume that
$\beta = 1$; otherwise, replace the function $f(x)$ with the
function $\beta^{n/2} f( \beta^{1/2} x )$. In this proof, $c, C,
C^{\prime}$ etc. stand for constants depending only on $A$ and
$\alpha$. We begin by proving (ii). Since $\sqrt{\EE |X|^2} =
\sqrt{n}$, Lemma \ref{lem_1010}(i) implies that
$$
 Prob \left \{ |X| \geq t \sqrt{n} \right \}  \leq 2 e^{-t/10}
\ \ \ \text{for all} \ t > 0.
$$
Therefore,
\begin{eqnarray}
\label{eq_342} \lefteqn{ \left| n - r^2 \right| \leq \EE \left|
|X|^2 - r^2 \right | = \int_0^{\infty} Prob \left \{ \left| |X|^2
- r^2 \right| \geq t \right \} dt } \\ & \leq & \int_0^{r^2} A
\exp \left( -\frac{\alpha t^2 n}{8 r^4} \right)dt +
\int_{r^2}^\infty \min \left \{  A e^{-\alpha n}, 2 \exp \left(
-\frac{\sqrt{t}}{10 \sqrt{ n}} \nonumber
 \right) \right \} dt \\
 & \leq & C \frac{r^2}{\sqrt{n}} + C^{\prime} n^3 A e^{-\alpha n}  + C^{\prime \prime}
 e^{-c n} < C \frac{r^2}{\sqrt{n}} + \tilde{C} e^{-\tilde{c} n},
\nonumber
\end{eqnarray}
provided that $n > C$. From (\ref{eq_342}) we deduce  (ii). To
prove (i), it is enough to consider the case where $\eps \geq
\frac{C}{\sqrt{n}}$. In this case, by (ii),
$$ Prob \left \{ \left| |X| - \sqrt{n} \right|
\geq \eps \sqrt{n} \right \}
\leq
Prob \left \{ \big| |X| - r \big|
\geq C^{\prime} \eps r \right \} $$
and (i) follows from (\ref{eq_855})
for the range $0 < \eps < 1 / C^{\prime}$. By adjusting the constants,
we establish (i) for the entire range $0 \leq \eps \leq 1$.
 \hfill $\square$

\begin{lemma} Let $n \geq 1$ be an integer, let $\beta > 0$, and let
$f: \RR^n \rightarrow [0, \infty)$ be a log-concave function that
is the density of a random vector with zero mean and with
covariance matrix that equals $\beta Id$. Then
$$ f(0) \geq e^{-n} \sup_{x \in \RR^n} f(x) \geq \left( \frac{c}{\sqrt{\beta}} \right)^n $$
where $0 < c < 1$ is a universal constant.
\label{zero_large_lem}
\end{lemma}

\emph{Proof:} The inequality $f(0) \geq e^{-n} \sup f$
is proved in  \cite[Theorem 4]{fradelizi}.
By our assumptions, $\int_{\RR^n} |x|^2 f(x) dx =
\beta n$. Markov's inequality entails
$$
 \int_{\sqrt{2 \beta n} D^n} f(x) dx \geq \frac{1}{2}.
$$ Therefore,
$$ \sup_{x \in \RR^n} f \geq
 \frac{1}{Vol(\sqrt{2 \beta n} D^n)}
 \int_{\sqrt{2 \beta n} D^n} f(x) dx  \geq (C \beta)^{-n/2}
\cdot \frac{1}{2}, $$ since $Vol (\sqrt{n} D^n) \leq \tilde{C}^n$
(see, e.g., \cite[page 11]{pisier_book}).  \hfill $\square$

\smallskip
\emph{Proof of Proposition \ref{cor_735}:}
Recall our assumption (\ref{eq_848_}) and our desired conclusion
(\ref{eq_309}) from the formulation of the proposition.
 We assume that $n$ is greater than
some large universal constant, since otherwise (\ref{eq_309})
 is obvious for an appropriate choice of constants $C, c > 0$.
 Denote $g = f * \gamma_{n,1}$, the convolution
of $f$ and $\gamma_{n,1}$.  Then $g$ is log-concave, and is the
density of a random vector with mean zero and covariance matrix $2
Id$. By Lemma \ref{zero_large_lem},
\begin{equation}
g(0) \geq \bar{c}^n. \label{zero_large}
\end{equation}
We set $C_0 = 25
\left(1 + \log 1 / \bar{c} \right) $ where $0 < \bar{c} < 1$ is the
constant from (\ref{zero_large}). Our assumption (\ref{eq_848_})
is precisely the basic requirement of Lemma \ref{lem_fourier}, for
$\alpha = C_0 / 5 \geq 5$. By the conclusion of that lemma,
\begin{equation}
\sup_{\theta \in S^{n-1}} g(t \theta) \leq e^{-5 n} g(0) +
\inf_{\theta \in S^{n-1}} g(t \theta) \ \ \ \text{for all} \ \ t
\geq 0, \label{eq_742_}
\end{equation}
since $e^{-C_0 n / 5} \leq e^{-5 n} g(0)$, according to the
definition of $C_0$ and (\ref{zero_large}). The function $g$ is
$C^{\infty}$-smooth, since $g = f * \gamma_{n,1}$ with
$\gamma_{n,1}$ being $C^{\infty}$-smooth. Additionally,
 since $0 < \int g < \infty$ then for some $A, B > 0$,
\begin{equation}
 g(x) \leq A e^{-B |x|} \ \ \ \text{for all} \ \ x \in \RR^n
\label{eq_408}
\end{equation}
(see, e.g., \cite[Lemma 2.1]{psitwo}). For $\theta \in S^{n-1}$
and $t \geq 0$, we write $g_{\theta}(t) = g(t \theta)$. Then
$g_{\theta}$ is log-concave, continuous on $[0, \infty)$,
$C^{\infty}$-smooth on $(0, \infty)$ and integrable on $[0,
\infty)$ by (\ref{eq_408}). In addition, $g_{\theta}(0) = g(0) >
0$  by (\ref{zero_large}). Fix $\theta_0 \in S^{n-1}$, and denote
$r_0 = t_n(g_{\theta_0})$. According to (\ref{eq_742_}), for any
$\theta \in S^{n-1}$ and $t \geq 0$,
$$
 |g_{\theta}(t) -  g_{\theta_0}(t)| \leq e^{-5 n} g(0) = e^{-5 n} \min \{ g_{\theta} (0), g_{\theta_0}(0) \}.
 $$
Thus the functions  $g_{\theta}$ and $g_{\theta_0}$ satisfy the
assumptions of Lemma \ref{lem_257}, for any $\theta \in S^{n-1}$.
By the conclusion of that lemma, for any $\theta \in S^{n-1}$,
$$
(1 - e^{-n}) r_0 \leq  t_n(g_{\theta}) \leq (1 +  e^{-n}) r_0,
$$
because $r_0 =  t_n(g_{\theta_0})$.
We deduce that for any
$10
e^{-n} \leq \eps \leq 1$ and $\theta \in S^{n-1}$,
\begin{equation}
 (1 + \eps) r_0 \geq \left(1 + \frac{\eps}{2} \right)
t_n(g_{\theta}) \ \ \ \text{and} \ \ \ (1 - \eps) r_0 \leq \left(1
- \frac{\eps}{2} \right) t_n(g_{\theta}).
\label{eq_303}
\end{equation}
For $0 \leq \eps \leq 1$ let $A_{\eps} = \{ x \in \RR^n ; \big |
|x| -r_0 \big| \leq \eps r_0 \}$. We will prove
that for all $0 \leq \eps \leq 1$,
\begin{equation}
\int_{A_{\eps}} g(x) dx \geq 1 - C e^{-c \eps^2 n}. \label{eq_313}
\end{equation}
Note that  (\ref{eq_313}) is obvious for $\eps < 10 e^{-n} \leq
\frac{10}{\sqrt{n}}$, since in this case $1 - C e^{-c \eps^2 n}
\leq 0$ for an appropriate choice of universal constants $c, C
> 0$. We still need to deal with the case
$10 e^{-n} \leq \eps \leq 1$. To that end, note that $g_{\theta}$
satisfies the requirements of Lemma \ref{lem_1153} for any $\theta
\in S^{n-1}$ by the discussion above. We will integrate in polar
coordinates and use (\ref{eq_303}) as well as Lemma
\ref{lem_1153}. This yields
\begin{eqnarray*}
 \lefteqn{ \int_{A_{\eps}} g(x) dx \geq
\int_{S^{n-1}} \int_{(1 - \eps/2) t_n(g_\theta)}^{(1 +
\eps/2) t_n(g_\theta)} t^{n-1} g_{\theta}(t) dt d\theta } \\
& \geq & \left(1 - C e^{-c \eps^2 n} \right) \int_{S^{n-1}}
\int_{0}^{\infty} t^{n-1} g_{\theta}(t) dt d\theta = 1 - C e^{-c \eps^2
n},
\end{eqnarray*}
since $\int_{\RR^n} g = 1$. This completes the proof of (\ref{eq_313}).

\smallskip
Let $X_1,X_2,...$ be a sequence of independent, real-valued,
standard gaussian random variables. By the classical central limit
theorem,
$$ Prob \left \{ \sum_{i=1}^m X_i^2 \leq m \right \} \stackrel{m \rightarrow \infty}{\longrightarrow} \frac{1}{2}. $$
Consequently, $1 / C^{\prime} \leq Prob \{ \sum_{i=1}^n X_i^2 \leq  n\} \leq 1
- 1 /C^{\prime} $ for some universal constant $C^{\prime} > 0$.
 Denote $X =
(X_1,...,X_n)$. Then $X$
is distributed according
to the density $\gamma_{n,1}$ in $\RR^n$. We record the bound just
mentioned: \begin{equation} \frac{1}{C^{\prime}} \leq Prob \{
|X|^2 \leq n \} \leq 1 - \frac{1}{C^{\prime}}. \label{bound_X}
\end{equation}
Let $Y$ be another random vector in $\RR^n$, independent of $X$,
that is distributed according to the density $f$. Since the
density of $X$ is an even function, then for any measurable sets
$I, J \subset [0, \infty)$ with $Prob \{ |X| \in I \} > 0$ and
$Prob \{ |Y| \in J \} > 0$,
\begin{equation} Prob \left \{ \langle
X, Y \rangle \geq 0 \ \text{\it given that} \  |X| \in I, |Y| \in
J \right \} = \frac{1}{2}. \label{eq_130}
\end{equation}
Additionally, the random vector $X + Y$ has $g$ as
its density, because $g = f * \gamma_{n,1}$. Therefore
(\ref{eq_313}) translates to
\begin{equation}
Prob \left \{ \big | |X + Y| - r_0 \big | > \eps r_0 \right \}
\leq C e^{-c \eps^2 n } \ \ \ \text{for all} \ \ 0 \leq \eps \leq
1. \label{eq_132}
\end{equation}
 Since
$X$ and $Y$ are independent, we conclude from (\ref{bound_X}),
(\ref{eq_130}) and (\ref{eq_132}) that for all $0 < \eps < 1$,
\begin{eqnarray}
\label{eq_136} \lefteqn{  Prob \left \{ |Y|^2 \geq r_0^2 (1 +
\eps)^2 - n \right \} } \\ & \leq & 2 C^{\prime} Prob \left \{
|Y|^2 \geq r_0^2 (1 + \nonumber
\eps)^2 - n, |X| \geq \sqrt{n}, \langle X, Y \rangle \geq 0 \right \}  \\
& \leq & 2 C^{\prime} Prob \left \{ |X + Y|^2 \geq r_0^2 (1 +
\eps)^2 \right \} \leq C \exp \left( -c \eps^2 n \right),
\nonumber
\end{eqnarray}
and similarly,
\begin{eqnarray}
\label{eq_141} \lefteqn{  Prob \left \{ |Y|^2 \leq r_0^2 (1 -
\eps)^2 - n \right \} } \\ & \leq & 2 C^{\prime} Prob \left \{
|Y|^2 \leq r_0^2 (1 - \nonumber
\eps)^2 - n, |X| \leq \sqrt{n}, \langle X, Y \rangle \leq 0 \right \}  \\
& \leq & 2 C^{\prime} Prob \left \{ |X + Y| \leq r_0 (1 - \eps)
\right \} \leq C \exp \left( -c \eps^2 n \right). \nonumber
\end{eqnarray}
Next, we estimate $r_0$. Recall that the density of $X+Y$ is
log-concave, $\EE(X + Y) = 0$ and $Cov(X + Y) = 2 Id$. We invoke
Lemma \ref{crude_lem}(ii), based on (\ref{eq_132}), and conclude
that $3 n / 2\leq r_0^2 \leq 3 n$, under the legitimate assumption
that $n
> C$. Denote $r = \sqrt{r_0^2 - n}$. Then $ \sqrt{n/2} \leq r
\leq \sqrt{2n}$ and
$$ r^2(1 + 10 \eps)^2 \geq r_0^2 (1 +
\eps)^2 - n,  \ \  \ \ \ r_0^2 (1 - \eps)^2 - n
\geq r^2 (1 - 10 \eps)^2,
$$
for $0 \leq \eps \leq 1/10$.
Therefore, (\ref{eq_136}) and (\ref{eq_141}) imply that for any $0
< \eps < \frac{1}{10}$,
$$
 Prob \left \{  r^2 (1 - 10 \eps)^2 \leq  |Y|^2  \leq r^2 (1 + 10
\eps)^2 \right \} \geq 1 - 2 C e^{-c \eps^2 n}.
$$
After adjusting the constants, we see that
\begin{equation}
\forall 0 \leq \eps \leq 1, \ \ \ \  Prob \left \{  \left| \, \frac{|Y|}{r} - 1 \, \right| \geq \eps
\right \} \leq C^{\prime} e^{-c^{\prime} \eps^2 n}.
\label{eq_346}
\end{equation}
Recall that $Y$ is distributed according to the
density $f$, which is an isotropic, log-concave function.
We may
thus apply Lemma \ref{crude_lem}(i), based on (\ref{eq_346}), and
conclude (\ref{eq_309}). The proposition is proved. \hfill
$\square$

\smallskip We proceed to discuss applications
of Proposition \ref{cor_735}. The following lemma is  usually
referred to as the Johnson-Lindenstrauss dimension reduction lemma
\cite{JL}. We refer, e.g., to \cite[Lemma 2.2]{gupta} for an
elementary proof.
 Recall that we denote by $Proj_E(x)$ the orthogonal
projection of $x$ onto $E$, whenever $x$ is a point in $\RR^n$ and
$E \subset \RR^n$ is a subspace.

\begin{lemma} Let $1 \leq k \leq n$  be integers, and let $E \in
G_{n,k}$ be a random $k$-dimensional subspace. Let $x \in \RR^n$
be a fixed vector. Then for all $0 \leq \eps \leq 1$,
\begin{equation}
Prob \left \{ \left | \, |Proj_E(x)| - \sqrt{\frac{k}{n}} |x| \,
\right | \geq \eps \sqrt{\frac{k}{n}} |x| \right \} \leq C e^{-c
\eps^2 k }
\label{eq_243}
\end{equation}
where $c, C > 0$ are universal constants.
\label{lem_find_reference}
\end{lemma}

\smallskip \emph{Proof of Theorem \ref{cor_202}:} We use the constant
$C_0 \geq 1$ from Proposition \ref{cor_735}, and the constant $c$
from Lemma \ref{lem_dvo}. Let $\ell = \lfloor \frac{c}{100 C_0}
\log n \rfloor$ and  fix $0 \leq \eps \leq 1 / 3$. We may assume
that $\ell \geq 1$; otherwise, $n$ is smaller than some universal
constant and the conclusion of the theorem is obvious. We assume
that $X$ is a random vector in $\RR^n$ whose density is an
isotropic, log-concave function to be denoted by $f$. Let $E \in
G_{n, \ell}$ be a fixed subspace that satisfies
\begin{equation}
\sup_{\theta \in S^{n-1} \cap E} M_{f}(\theta, t) \leq e^{-C_0
\ell} + \inf_{\theta \in S^{n-1} \cap E} M_{f}(\theta, t) \ \ \
\text{for all} \ \ t \in \RR. \label{eq_609_}
\end{equation}
Denote $g = \pi_E(f)$.  Then (\ref{eq_609_}) translates, with the
help of (\ref{proj_ok}) from Section \ref{section2}, to
\begin{equation}
\sup_{\theta \in S^{n-1} \cap E} M_{g}(\theta, t) \leq e^{-C_0
\ell} + \inf_{\theta \in S^{n-1} \cap E} M_{g}(\theta, t) \ \ \
\text{for all} \ \ t \in \RR. \label{eq_609}
\end{equation}
The function $g$ is an isotropic, log-concave function, and it is
the density of $Proj_E(X)$. We invoke Proposition \ref{cor_735},
for $\ell$ and $g$, based on (\ref{eq_609}). By the conclusion of
that proposition,
\begin{equation}
Prob \left \{ \left| \, \frac{|Proj_E(X)|}{\sqrt{\ell}} - 1  \,
\right| \geq \eps \right \} \leq C^{\prime} e^{-c^{\prime} \eps^2
\ell}, \label{eq_625}
\end{equation}
under the assumption that the subspace $E$ satisfies
(\ref{eq_609_}). Suppose that $F \in G_{n,\ell}$ is a random
$\ell$-dimensional subspace in $\RR^n$, independent of  $X$.
Recall our choice of the integer $\ell$. According to Lemma
\ref{lem_dvo}, with probability greater than $1 - e^{-c
n^{0.99}}$, the subspace $E = F$ satisfies (\ref{eq_609_}). We
conclude from (\ref{eq_625}) that
$$
Prob \left \{ \left| \, \frac{|Proj_F(X)|}{\sqrt{\ell}} - 1  \,
\right| \geq \eps \right \} \leq C^{\prime} e^{-c^{\prime} \eps^2
\ell} + e^{-c n^{0.99}} \leq \tilde{C} e^{-\tilde{c} \eps^2 \ell},
$$
where the last inequality holds as $\ell \leq \log n$
and $0 \leq \eps \leq 1/3$. Since $X$
and $F$ are independent, then by Lemma \ref{lem_find_reference},
$$ Prob \left \{ \left | \, |Proj_F(X)| -
\sqrt{\frac{\ell}{n}} |X| \, \right | \geq \eps
\sqrt{\frac{\ell}{n}} |X| \right \} \leq \hat{C} e^{-\hat{c}
\eps^2 \ell }. $$ To summarize, with probability greater than $1 -
\bar{C} e^{-\bar{c} \eps^2 \ell }$ we have
\begin{enumerate}
\item[(i)] $\displaystyle \ \ \ (1 - \eps) \sqrt{\ell} \leq |Proj_F(X)| \leq  (1 +
\eps) \sqrt{\ell}$, and  also
\item[(ii)] $\displaystyle  \ \ \ (1 + \eps)^{-1}
\sqrt{\frac{n}{\ell}} |Proj_F(X)| \leq |X| \leq (1 - \eps)^{-1}
\sqrt{\frac{n}{\ell}} |Proj_F(X)|$.
\end{enumerate}
Hence,
\begin{equation}
  Prob \left \{ \frac{1 - \eps}{1+\eps} \leq
\frac{|X|}{\sqrt{n}} \leq \frac{1 + \eps}{1-\eps} \right \} \geq 1
- \bar{C} e^{-\bar{c} \eps^2 \ell}. \label{eq_841} \end{equation}
Note that $\frac{1+\eps}{1-\eps} \leq 1 + 3 \eps$ and  $1 - 3 \eps
\leq \frac{1-\eps}{1+\eps}$, and recall that $0 \leq \eps \leq
\frac{1}{3}$ was arbitrary, and that $\ell = \lfloor \frac{c}{100
C_0} \log n \rfloor$. By adjusting the constants, we deduce from
(\ref{eq_841}) that the inequality in the conclusion of the
theorem is valid for all $0 \leq \eps \leq 1$. The theorem is thus
proved. \hfill $\square$

\smallskip The following lemma may be proved via a straightforward
computation. Nevertheless, we will present a shorter, indirect
proof that is based on properties of the heat kernel, an idea we
borrow from \cite[Theorem 3.1]{brehm_voigt}.

\begin{lemma} Let $n \geq 1$ be an integer
 and let $\alpha, \beta > 0$. Then,
\begin{equation}
\int_{\RR^n} \left|\gamma_{n,\alpha}(x) - \gamma_{n, \beta}(x)
\right| dx \leq C \sqrt{n} \left| \frac{\beta}{\alpha} - 1 \right|,
\label{eq_1021}
\end{equation}
where $C > 0$ is a universal constant. \label{lem_929}
\end{lemma}

\emph{Proof:} The integral on the left-hand side of
(\ref{eq_1021}) is never larger than $2$. Consequently, the lemma
is obvious when $\frac{\beta}{\alpha} > 2$ or when
$\frac{\beta}{\alpha} < \frac{1}{2}$, and hence we may assume that
$\frac{1}{2} \alpha \leq \beta \leq 2 \alpha$. Moreover, in this
case both the left-hand side and the right-hand side of
(\ref{eq_1021}) are actually symmetric in $\alpha$ and $\beta$
 up
to a factor of at most $2$. Therefore, we may assume that $\alpha
< \beta \leq 2 \alpha$ (the case $\beta = \alpha$ is obvious). For
$t
> 0$ and for a measurable function $f: \RR^n \rightarrow \RR$, we define
$$ (P_t f)(x) = \frac{1}{(4 \pi t)^{n/2}} \int_{\RR^n}
e^{-\frac{|x-y|^2}{4 t}} f(y) dy \ \ \ \ \ (x \in \RR^n) $$
whenever the integral converges. Then $(P_t)_{t > 0}$ is the heat
semigroup on $\RR^n$. We will make use of the following estimate:
For any smooth, integrable function $f: \RR^n \rightarrow \RR$ and
any $t >0$,
\begin{equation}
\int_{\RR^n} \left| (P_t f)(x) - f(x) \right| dx \leq 2 \sqrt{t}
\int_{\RR^n} |\nabla f(x)| dx. \label{eq_1032_}
\end{equation}
An elegant proof of the inequality (\ref{eq_1032_}), in a much
more general setting, is given by Ledoux \cite[Section
5]{ledoux}. It is straightforward to verify that
\begin{eqnarray*}
\lefteqn{ \int_{\RR^n} |\nabla \gamma_{n, \alpha}(x)| dx
= \frac{1}{(2 \pi \alpha)^{n/2}} \int_{\RR^n} \frac{|x|}{\alpha}
e^{-\frac{|x|^2}{2 \alpha}} dx } \\ & \leq & \frac{1}{\alpha} \left(
\frac{1}{(2 \pi \alpha)^{n/2}} \int_{\RR^n} |x|^2
e^{-\frac{|x|^2}{2 \alpha}} dx \right)^{1/2} =
\sqrt{\frac{n}{\alpha}}.  \phantom{aaaaaaaaaa}
\end{eqnarray*}
Consequently, (\ref{eq_1032_}) implies  that
\begin{equation}
\int_{\RR^n} \left| P_{\frac{\beta - \alpha}{2}}
\left( \gamma_{n, \alpha} \right)  (x)  - \gamma_{n, \alpha}(x)
\right| dx \leq 2 \sqrt{ \frac{\beta - \alpha}{2} }
\sqrt{\frac{n}{\alpha}}.
\label{eq_534}
\end{equation}
It is well-known and easy to prove that $ \gamma_{n, \beta} =
P_{\frac{\beta - \alpha}{2}} \left( \gamma_{n, \alpha} \right)$.
Since $\alpha < \beta \leq 2 \alpha$, then (\ref{eq_534}) implies
(\ref{eq_1021}). The lemma is proved. \hfill $\square$

\smallskip We are now able to prove Theorem \ref{thm_uncond} by combining
the classical Berry-Esseen bound with Theorem \ref{cor_202}.

\smallskip
\emph{Proof of Theorem \ref{thm_uncond}:} We may assume that $n$
exceeds a given universal constant. Let $f$ and $X$ be as in the
assumptions of Theorem \ref{thm_uncond}. According to Theorem
\ref{cor_202},
\begin{equation}
Prob \left \{ \left| \,  \frac{|X|}{\sqrt{n}} - 1 \, \right| \geq
\eps \right \}  \leq C n^{-c \eps^2} \ \ \ \text{for all} \ \ 0
\leq \eps \leq 1. \label{eq_115_}
 \end{equation}
The case $\eps = \sqrt{2}-1$ in (\ref{eq_115_}) shows that
$\delta_0 := Prob \left \{ |X| \geq \sqrt{2 n} \right \} \leq C
n^{-c/4} \leq n^{-c/10}$, under the legitimate assumption that $n$
exceeds a certain universal constant. By (\ref{eq_115_}) and by
Lemma \ref{lem_1010}(ii),
\begin{eqnarray}
\label{eq_115__} \lefteqn{ \EE \left| \,  \frac{|X|^2}{n} - 1 \,
\right| = \int_0^{\infty}  Prob \left \{ \left| \, \frac{|X|^2}{n}
- 1 \, \right| \geq t \right \} dt } \\ & \leq & \int_0^1
C^{\prime} n^{-c^{\prime} t^2} dt + \int_1^{\infty} (1 - \delta_0)
 \left( \frac{\delta_0}{1 - \delta_0}
\right)^{(\sqrt{\frac{1+t}{2}} + 1)/2} dt \leq \frac{C^{\prime
\prime}}{\sqrt{\log n}}. \nonumber
 \end{eqnarray}
Let $\delta_1,...,\delta_n$ be independent Bernoulli random
variables, that are also independent of $X$, such that $Prob \{
\delta_i = 1 \} = Prob \{ \delta_i = -1 \} = 1/2$ for $i=1,...,n$.
 For $t \in \RR$ and $x = (x_1,...,x_n) \in \RR^n$
denote
$$ P(x;t) = Prob \left \{ \frac{\sum_{i=1}^n \delta_i x_i}{\sqrt{n}}
\leq t \right \}. $$ We write $$ \Phi_{\sigma^2}(t) =
\frac{1}{\sqrt{2 \pi \sigma}} \int_{-\infty}^t \exp \left( -
\frac{t^2}{2 \sigma^2} \right) dt
$$ for $\sigma > 0$ and $t \in \RR$. By the Berry-Esseen bound
(see, e.g., \cite[Section XVI.5]{feller} or \cite[Section
2.1.30]{stroock}), for any $x \in \RR^n$,
\begin{equation}
\sup_{t \in \RR} \left| \, P(x; t) \, - \, \Phi_{|x|^2 / n}(t) \,
\right| \leq C \frac{\sum_{i=1}^n |x_i|^3}{|x|^3}, \label{eq_1230}
\end{equation}
where $C > 0$ is a universal constant. Since $f$ is unconditional,
the random variable $ \left( \sum_{i=1}^n X_i \right) / \sqrt{n}$
has the same law of distribution as the random variable $\left(
\sum_{i=1}^n \delta_i X_i \right) / \sqrt{n}$. For $t \in \RR$ we
set
$$ P(t) = Prob \left \{ \frac{\sum_{i=1}^n X_i}{\sqrt{n}}
\leq t \right \} = Prob \left \{ \frac{\sum_{i=1}^n \delta_i
X_i}{\sqrt{n}} \leq t \right \} . $$ We denote the expectation
over the random variable $X$ by $\EE_X$. Then $ P(t) = \EE_X P(X;
t) $ by the complete probability formula.  For $i=1,...,n$, the
random variable $X_i$ has mean zero, variance one, and its density
is a log-concave function. Consequently, $\EE |X_i|^2 = 1$, and by
Lemma \ref{lem_1010}(i), for any $1 \leq i \leq n$,
$$
 Prob \left \{ |X_i| \geq 20 \log n \right \} \leq 2 e^{-2 \log n} =
\frac{2}{n^2}. $$ Therefore, with probability greater than $1 -
\frac{2}{n}$ of selecting $X$,
\begin{equation}
|X_i| \leq 20 \log n \ \ \ \text{for all} \ \ 1 \leq i \leq n.
\label{eq_109}
\end{equation}
Fix $t \in \RR$.  We substitute into (\ref{eq_1230}) the
information from (\ref{eq_109}), and from the case $\eps = 1/2$ in
(\ref{eq_115_}). We see that with probability greater than $1 - C
n^{-c/4} - \frac{2}{n}$ of selecting $X$,
$$
 \left| \, P(X;t) \, - \, \Phi_{\frac{|X|^2}{
n}}(t) \, \right| \leq   C \frac{\sum_{i=1}^n |X_i|^3}{|X|^3} \leq
C^{\prime} \frac{(\log n)^3}{\sqrt{n}}.
$$
Since always $0 \leq P(X ;t) \leq 1$ and $0 \leq \Phi_1(t) \leq
1$, we conclude that
\begin{equation} \EE_X \left| \, P(X;t) \, - \,
\Phi_{\frac{|X|^2}{ n}}(t) \, \right| \leq C^{\prime} \frac{(\log
n)^3}{\sqrt{n}} + 2 C n^{-c/4} + \frac{2}{n} <
\frac{C^{\prime}}{n^{c^{\prime}}}. \label{eq_117}
\end{equation}
According  to Lemma \ref{lem_929}, for any $x \in \RR^n$,
$$
\left| \, \Phi_{\frac{|x|^2}{n}}(t) \, - \, \Phi_1(t) \, \right|
\leq \int_{-\infty}^{\infty} \left| \gamma_{1,\frac{|x|^2}{n}}(s)
- \gamma_{1,1}(s) \right| ds \leq
 \hat{C} \left| \, \frac{|x|^2}{n} - 1
\, \right |,
$$
and therefore by (\ref{eq_115__})
\begin{equation}
 \EE_X \left| \, \Phi_{\frac{|X|^2}{n}}(t) \, - \, \Phi_1(t) \,
\right| \leq  \hat{C} \EE_X \left| \, \frac{|X|^2}{n} - 1 \,
\right | \leq \frac{C}{\sqrt{\log n}}. \label{eq_609__}
\end{equation}
Recall that $ P(t) = \EE_X P(X; t) $ and that $t$ is an arbitrary
real number. We apply Jensen's inequality, and then combine
(\ref{eq_117}) and (\ref{eq_609__}) to obtain
\begin{equation}
\forall t \in \RR, \ \ \ \left| \, P(t) \, - \, \Phi_{1}(t) \,
\right| \leq \EE_X \left| \, P(X ; t) \, - \, \Phi_{1}(t) \,
\right| \leq \frac{C}{\sqrt{\log n}}. \label{eq_200}
\end{equation}
The random variable  $(X_1 + ... + X_n) / \sqrt{n}$ has mean zero,
variance one and a log-concave density. Its cumulative
distribution function $P(t) = Prob \left \{ (X_1 + ... + X_n) /
\sqrt{n} \leq t \right \}$ satisfies (\ref{eq_200}). Therefore, we
may invoke \cite[Theorem 3.3]{BHVV}, and conclude from
(\ref{eq_200}) that
$$ d_{TV} \left(  \, \frac{X_1 + ... + X_n}{\sqrt{n}} \, ,
\, Z \, \right) \leq \check{C} \left(
\frac{C \log \frac{C}{\sqrt{\log n}}}{\sqrt{\log n}} \right)^{1/2}
 = \check{C}
\frac{\sqrt{\log \log n}}{(\log n)^{1/4}}, $$
 where $Z \sim N(0,1)$ is a standard gaussian random variable.
The theorem follows, with $\eps_n \leq C (\log \log (n+2))^{1/2} /
(\log (n+1) )^{1/4}$. \hfill $\square$

\smallskip \emph{Remarks.} \begin{enumerate}
\item Suppose that $f$ is a log-concave density in high dimension
that is isotropic
and unconditional. In Theorem \ref{thm_uncond}, we were able to
describe an explicit one-dimensional marginal of $f$ that is
approximately normal. It seems possible to identify some
multi-dimensional subspaces $E \subset \RR^n$, spanned by specific
sign-vectors, such that $\pi_E(f)$ is guaranteed to be
almost-gaussian. We did not pursue this direction.

\vspace{5pt}
\item
Under the assumptions of Theorem \ref{thm_uncond}, we proved that
$\langle X, \theta \rangle$ is approximately gaussian
 when $\theta = (1,...,1) / \sqrt{n}$.
A  straightforward adaptation of
the proof of Theorem \ref{thm_uncond} shows
that
 $\langle X, \theta \rangle$ is approximately gaussian under the weaker assumption that
$|\theta_1|,...,|\theta_n|$ are rather small (as in Lindeberg's
condition).

\vspace{5pt}
\item Theorem \ref{thm_basic}, with a worse bound for
$\eps_n$, follows by combining Theorem \ref{cor_202} with
 the methods in \cite{ABP}, and then applying  \cite[Theorem 3.3]{BHVV}.
We will deduce Theorem \ref{thm_basic} from the stronger Theorem
\ref{thm_multi} in the next section.
\end{enumerate}

\setcounter{equation}{0}
\section{Multi-dimensional marginals} \label{section5}
 The next few pages
are devoted to the proof of the following lemma.

\begin{lemma}
Let $n \geq 2$ be an integer, let $\alpha \geq 10$, and let $f:
\RR^n \rightarrow [0, \infty)$ be an isotropic, log-concave
function. Denote $g = f * \gamma_{n,n^{-30 \alpha}}$. Then,
$$ \int_{\RR^n} |g(x) - f(x)| dx \leq \frac{C}{n^{\alpha/10}}, $$
where $C > 0$ is a universal constant.
\label{lem_915}
\end{lemma}

We begin with an addendum to Lemma \ref{lem_1153}. Rather than
appealing to the Laplace asymptotic method once again, we will
base our proof on an elegant observation by Bobkov regarding
one-dimensional log-concave functions.

\begin{lemma} Let $n \geq 2$ be an integer, let $\alpha \geq 5$
and let $f: [0, \infty) \rightarrow [0, \infty)$ be a log-concave
function with $\int f < \infty$. Denote $t_0 = \sup \{ t > 0 ;
f(t) \geq e^{-\alpha n} f(0) \}$. Then,
\begin{equation}
\int_0^{t_0} t^{n-1} f(t) dt \geq \left(1 - e^{-\alpha n / 8}
\right) \int_{0}^{\infty} t^{n-1} f(t) dt. \label{eq_329}
\end{equation} \label{lem_1153_}
\end{lemma}

\emph{Proof:} If $\int f = 0$ then $f \equiv 0$ almost everywhere
and (\ref{eq_329}) is trivial. Thus, we may suppose that $\int f >
0$. Moreover, we may assume that $f$ is continuous on $[0,
\infty)$ and $C^2$-smooth on $(0, \infty)$, by approximation (for
example, convolve $f$ with $\gamma_{1, \eps}$ on $\RR$, restrict
the result to $[0, \infty)$, and let $\eps$ tend to zero). Since
$0 < \int f < \infty$ then $f$ decays exponentially fast at
infinity, and $0 < \int_0^{\infty} t^{n-1} f(t) dt < \infty$.
Multiplying $f$ by a positive constant, we may assume that
$\int_0^{\infty} t^{n-1} f(t) dt = 1$.

\smallskip For $t > 0$, denote,
$$ \phi(t) = t^{n-1} f(t) \ \ \ \text{and} \ \ \ \Phi(t) = \int_0^{t} \phi(s)
ds. $$ Then $\phi$ is a log-concave function with $\int \phi = 1$.
Recall the definition of $t_n(f)$, that is, (\ref{eq_106}) from
Section \ref{section4}. According to that definition,
$\phi^{\prime}(t_n(f)) = 0$. Denote $M = f(t_n(f)) > 0$. Then $M
\geq e^{-(n-1)} f(0)$ by (\ref{eq_1127}) from Section
\ref{section4}, and hence
$$ t_0 \geq t_1 := \sup \left \{ t > 0; f(t) \geq e^{-(\alpha - 1) (n-1)} M \right \}, $$
where $t_0$ is defined in the formulation of the lemma. Since $M
> 0$ and since $f$ is continuous and vanishes at infinity, the
number $t_1$ is finite, greater than $t_n(f)$, and satisfies
$f(t_1) = e^{-(\alpha -1) (n-1)} M$. From (\ref{eq_1127_2}) of
Section \ref{section4} we see that $ t_1 \leq \alpha t_n(f)$.
Therefore,
\begin{eqnarray*}
\lefteqn{ \phi(t_1) = \phi(t_n(f)) \cdot \left( \frac{t_1}{t_n(f)}
\right)^{n-1} \cdot \frac{f(t_1)}{M} } \\ & \leq & \phi(t_n(f))
\cdot \alpha^{n-1} \cdot e^{-(\alpha - 1)(n-1)} \leq \phi(t_n(f))
e^{-\alpha n / 8 } = e^{-\alpha n / 8} \cdot \max \phi,
\end{eqnarray*} where $\phi(t_n(f)) = \max \phi$ because $\phi$ is
log-concave, $\phi(t_n(f)) > 0$ and $\phi^{\prime}(t_n(f)) = 0$.
Let $\Phi^{-1} : (0, 1) \rightarrow (0, \infty)$ stand for the
inverse function to $\Phi$. A useful fact we learned from Bobkov's
work \cite[Lemma 3.2]{bobkov} is that the function $\psi(t) =
\phi(\Phi^{-1}(t))$ is concave on $(0,1)$. (To see this,
differentiate $\psi$ twice, and use the inequality $(\log
\phi)^{\prime \prime} \leq 0$.)

Since $\phi$ attains its maximum at $t_n(f)$, then $\psi$ attains
its maximum at $\Phi(t_n(f))$. The function $\psi$ is non-negative
and concave on $(0,1)$, hence for $t \geq \Phi(t_n(f))$ and $0 <
\eps < 1$,
$$ \psi(t) \leq \eps \cdot \max \psi  \ \ \  \Rightarrow \ \ \ t \geq 1 - \eps.
$$
Equivalently, for $s \geq t_n(f)$ and $0 < \eps < 1$, the
inequality $\phi(s) \leq \eps \cdot \max \psi = \eps \cdot \max
\phi $ implies the bound $\Phi(s) \geq 1 - \eps$. We have shown
that $t_1 \geq t_n(f)$ satisfies $\phi(t_1) \leq e^{-\alpha n / 8}
\max \phi$, and hence we conclude that $\Phi(t_1) \geq 1 -
e^{-\alpha n / 8}$. Recalling that $t_0 \geq t_1$, the lemma
follows. \hfill $\square$

\begin{corollary} Let $n \geq 2$ be an integer, let $\alpha
\geq 5$, and let
$f: \RR^n \rightarrow [0, \infty)$ be a
log-concave function with $\int f = 1$.
Denote $K = \{ x \in \RR^n ; f(x) \geq e^{-\alpha n} f(0) \}$.
Then,
$$ \int_{K} f(x) dx \geq 1 - e^{-\alpha n / 8}. $$ \label{cor_302}
\end{corollary}

\emph{Proof:} For $\theta \in S^{n-1}$ set
$$ I(\theta) = \{ t \geq 0 ; f(t \theta) \geq e^{-\alpha n} f(0)
\} = \{ t \geq 0 ; t \theta \in K \}. $$ By log-concavity,
$I(\theta)$ is a (possibly infinite) interval in $[0, \infty)$
containing zero. For $t \geq 0$ and $\theta \in S^{n-1}$ we denote
$f_{\theta}(t) = f(t \theta)$. Then $f_{\theta}$ is log-concave.
Since $\int f =1$, then, e.g., by \cite[Lemma 2.1]{psitwo} we know
that $f$ decays exponentially fast at infinity and $\int
f_{\theta} < \infty$. Next, we integrate in polar coordinates and
use Lemma \ref{lem_1153_}. This yields
\begin{eqnarray*}
\lefteqn{ \int_{K} f(x) dx = \int_{S^{n-1}} \int_0^{\sup
I(\theta)} t^{n-1} f_{\theta}(t)  dt d\theta } \\ & \geq &
 \left(1 - e^{-\alpha n/8} \right) \int_{S^{n-1}} \int_0^\infty t^{n-1} f_{\theta}(t)
dt d \theta = 1 - e^{-\alpha n/8}.
\end{eqnarray*} \hfill $\square$

\begin{lemma} Let $n \geq 1$ be an integer and
let $X$ be a random vector in $\RR^n$ with an isotropic,
log-concave density. Suppose that $K \subset \RR^n$ is convex with
 $Prob \{ X \in K \}
\geq \frac{9}{10}$. Then,
$$ \frac{1}{10} D^n \subset K. $$
\label{lem_lovasz}
\end{lemma}

\emph{Proof:} Assume the contrary. Since $K$ is convex, then there
exists $\theta \in S^{n-1}$ such that $ K \subset \left \{ x \in
\RR^n ; \langle x, \theta \rangle < 1 / 10 \right \}$. Hence,
\begin{equation}
 Prob \left \{
 \langle X, \theta \rangle \leq \frac{1}{10} \right
\}   \geq Prob \left \{ X \in K \right \} \geq \frac{9}{10}.
\label{eq_1018}
\end{equation}
 Denote $E = \RR \theta$, the one-dimensional line spanned by $\theta$, and let $g =
\pi_E(f)$. Then $g$ is log-concave and isotropic, hence $\sup g
\leq 1$ by (\ref{eq_547}) of Section \ref{section2}. Since $g$ is
the density of the random variable $\langle X, \theta \rangle$ and
$\sup g \leq 1$, then
\begin{equation}
 Prob \left \{ 0 \leq \langle X, \theta \rangle \leq \frac{1}{10}
 \right \} = \int_0^{1/10} g(t) dt \leq \frac{1}{10}.
 \label{eq_956}
 \end{equation}
An appeal to \cite[Lemma 3.3]{bobkov} -- a result that essentially
goes back to Gr\"unbaum and Hammer \cite{grunbaum} -- shows that
\begin{equation}
 Prob \{ \langle X, \theta \rangle < 0 \} \leq 1 -
\frac{1}{e} < \frac{4}{5}. \label{eq_957}
\end{equation}
After adding (\ref{eq_957}) to (\ref{eq_956}), we arrive at a
contradiction to (\ref{eq_1018}). This completes the proof. \hfill
$\square$

\smallskip For two sets $A, B \subset \RR^n$ we
write $A + B = \{ x + y; x \in A, y \in B \}$ and $A - B = \{ x -
y; x \in A, y \in B \}$ to denote their Minkowski sum and
difference.

\begin{lemma} Let $n \geq 2$ be an integer,
let $\alpha \geq 10$, and let $f: \RR^n \rightarrow [0, \infty)$
be an isotropic,  log-concave function. Consider the sets $K_0 =
\{ x \in \RR^n; f(x) \geq e^{-\alpha n} f(0) \}$ and $K = \{ x \in
\RR^n ; \exists y \not \in K_0 , |x - y| \leq n^{-3 \alpha} \}$.
Then,
$$
 \int_K f(x) dx \leq \frac{C}{n^{\alpha}} $$
where $C > 0$ is a universal constant.
\label{lem_1208}
\end{lemma}

\emph{Proof:} Let $\mu$ be the probability measure on $\RR^n$
whose density is $f$. By Corollary \ref{cor_302},
\begin{equation}
\mu( K_0 ) = \int_{K_0} f(x) dx \geq 1 - e^{-\alpha n / 8} \geq
\frac{9}{10}. \label{eq_304}
\end{equation}
The set $K_0$ is convex, since $f$ is log-concave. According to
(\ref{eq_304}) and Lemma \ref{lem_lovasz},
\begin{equation}
\frac{1}{10}D^n \subset K_0. \label{eq_159}
\end{equation}
By the definition, $K = (\RR^n \setminus K_0) +  n^{-3 \alpha }
D^n$. Since $D^n \subset -10 K_0$, then
\begin{equation} K
\subset (\RR^n \setminus K_0) - 10 n^{-3 \alpha } K_0  \subset
\RR^n \setminus \left(1 - n^{-2 \alpha } \right) K_0,
\label{eq_944} \end{equation} because  $K_0$ is convex and $10
n^{-3 \alpha} \leq n^{-2 \alpha}$. We use (\ref{eq_159}) and Lemma
\ref{zero_large_lem} for $\beta = 1$. This implies the estimate
\begin{equation}
\mu \left( \frac{D^n}{20} \right) =
 \int_{\frac{D^n}{20}} f(x) dx
\geq e^{-\alpha n} f(0) \cdot Vol \left( \frac{D^n}{20} \right)
\geq
 \left( \frac{c^{\prime} e^{-\alpha} }{\sqrt{n}} \right)^n, \label{eq_218}
\end{equation}
where we also used the standard estimate $Vol(D^n) \geq \left( c /
\sqrt{n} \right)^n$. The inclusion (\ref{eq_159}) and the
convexity of $K_0$ entail that
$$
\left( 2 n^{-2 \alpha} \right) \frac{D^n}{20} + \left(1 - 2 n^{-2
\alpha} \right) K_0 \subset \left(1 - n^{-2 \alpha} \right) K_0.
$$
Therefore, according to the Pr\'ekopa-Leindler inequality,
\begin{equation}
\mu \left( \left(1 - n^{-2 \alpha } \right) K_0 \right) \geq \mu
\left(  \frac{D^n}{20} \right)^{ 2 n^{-2 \alpha} } \cdot  \mu
\left( K_0 \right)^{1 - 2 n^{-2 \alpha}}. \label{eq_1015}
\end{equation}
We combine (\ref{eq_944}), (\ref{eq_1015}), (\ref{eq_218}) and
 (\ref{eq_304}) to obtain
\begin{eqnarray*}
\lefteqn{ \mu(K) \leq \mu \left( \RR^n \setminus \left(1 - n^{-2
\alpha } \right) K_0 \right) = 1 - \mu \left( \left(1 - n^{-2
\alpha } \right) K_0 \right) } \\ & \leq & 1 - \left( \left(
\frac{c^{\prime} e^{-\alpha} }{\sqrt{n}} \right)^n \right)^{ 2
n^{-2 \alpha} } \cdot \left( 1 - e^{-\alpha n / 8} \right)^{1 - 2
n^{-2 \alpha}} \leq  \frac{C^{\prime}}{n^{\alpha}},
\end{eqnarray*}
for some universal constant $C^{\prime} > 0$ (the verification of
the last inequality is elementary and routine). The lemma is thus
proved. \hfill $\square$

\smallskip
\emph{Proof of Lemma \ref{lem_915}:} By approximation, we may
assume that $f$ is continuously differentiable. Denote $\psi =
\log f$ (with $\psi = -\infty$ when $f = 0$). Then $\psi$ is a
concave function. Consider the sets $K_0 = \{ x \in \RR^n; f(x)
\geq e^{-\alpha n} f(0) \}$ and $K = \{ x \in \RR^n ; \exists y
\not \in K_0 , |x - y| < n^{-4 \alpha} \}$.
 The first step of the proof is to show that
\begin{equation}
\{ x \in K_0 ; |\nabla \psi(x)| > n^{5 \alpha} \} \subset
 K. \label{eq_317}
\end{equation}
Note that $f(0) > 0$ by \cite[Theorem 4]{fradelizi},
and hence $f(x) > 0$ for all $x \in K_0$. Consequently,
$\psi$ is finite on $K_0$, and $\nabla \psi$ is well-defined
on $K_0$. In order to prove (\ref{eq_317}),
let us pick $x \in K_0$ such that $|\nabla \psi(x)| > n^{5
\alpha}$. Set $\theta = \nabla \vphi(x) / |\nabla
\vphi(x)|$. To prove (\ref{eq_317}), it suffices
to show that
$$
x - n^{-4 \alpha} \theta \not \in K_0,
$$
by the definition of $K$. According to the definition of $K_0$, it is enough
to prove that
\begin{equation}
f \left( x - n^{-4 \alpha} \theta \right) < e^{-\alpha n} f(0).
\label{eq_551}
\end{equation}
We thus focus on proving (\ref{eq_551}). We may assume that $f(x -
n^{-4 \alpha } \theta)
> 0$ since otherwise (\ref{eq_551}) holds trivially.
By concavity, $\vphi(t) := \psi(x + t \theta) = \log f(x + t
\theta)$ is finite for $ -n^{-4 \alpha } \leq t \leq 0$, and
$$ \vphi^{\prime}(0) = \langle \nabla \psi(x), \theta \rangle
= |\nabla \psi(x)| > n^{5 \alpha}. $$ Since $\vphi$ is concave,
then $\vphi^{\prime}$ is non-increasing. Consequently,
$\vphi^{\prime}(t) > n^{5 \alpha}$ for $-n^{-4 \alpha} \leq t \leq
0$. Hence,
\begin{equation}
 \vphi(0)  - \vphi(-n^{-4 \alpha}) >
n^{5 \alpha} \cdot n^{-4 \alpha} = n^{\alpha} \geq \alpha n  + 1,
\label{eq_256}
\end{equation}
as $\alpha \geq 10$ and $n \geq 2$. Recall that $f(0) \geq e^{-n}
f(x)$ by \cite[Theorem 4]{fradelizi} and that $f(x + t \theta) =
e^{\vphi(t)}$. We conclude from (\ref{eq_256}) that $f(0) \geq
e^{-n} f(x) > e^{\alpha n} f(x - n^{-4 \alpha} \theta)$, and
(\ref{eq_551}) is proved. This completes the proof of
(\ref{eq_317}).

\smallskip
For $x \in \RR^n$ and $\delta > 0$ denote $B(x, \delta) = \{ y \in
\RR^n ; |y-x| \leq \delta \}$. Fix $x \in K_0$ such that $B(x,
n^{-3 \alpha }) \subset K_0$. Then for any $y \in B(x, n^{-10
\alpha})$ we have $y \not \in K$ and hence $|\nabla \psi(y)| \leq
n^{5 \alpha }$, by (\ref{eq_317}). Consequently,
$$
|\psi(y) - \psi(x)| \leq n^{5 \alpha} |x - y| \leq n^{-5 \alpha} \
\ \ \text{for all} \ \ y \in B(x, n^{-10 \alpha}). $$
Recalling that $f = e^{\psi}$, we obtain
\begin{equation}
|f(y) - f(x)| \leq 2 n^{-5 \alpha} f(x)  \ \ \ \text{for all} \ \
y \in B(x, n^{-10 \alpha }). \label{eq_659}
\end{equation}
We will also make use of the crude estimate
\begin{equation}
\int_{ \RR^n \setminus B(0, n^{-10 \alpha})} \gamma_{n,n^{-30
\alpha}}(x) dx  \leq 2 \exp (-n^{4\alpha} / 10)  \leq e^{-20
\alpha n}, \label{eq_659_}
\end{equation}
that follows, for example, from Lemma \ref{lem_1010}(i) as
$\sqrt{\int_{\RR^n} |x|^2 \gamma_{n,n^{-30 \alpha}}(x) dx} =
n^{1/2-15 \alpha}$. According to \cite[Theorem 4]{fradelizi},
\begin{equation}
\sup f \leq e^n f(0) \leq e^{(\alpha + 1)n} f(x), \label{eqn_521}
\end{equation}
since $x \in K_0$. Recall  that $g = f * \gamma_{n, n^{-30
\alpha}}$. We use (\ref{eq_659}), (\ref{eq_659_}) and
(\ref{eqn_521}) to conclude that
\begin{eqnarray}
\label{eq_738} \lefteqn{ |g(x) - f(x)| \leq \int_{\RR^n}
\gamma_{n,n^{-30 \alpha}}(x-y) \left| f(y) - f(x) \right| dy } \\
& \leq & 2 n^{-5 \alpha } f(x) + 2 \sup f  \cdot \int_{\RR^n
\setminus B(x, n^{-10 \alpha})} \gamma_{n,n^{-30 \alpha}} (x-y) dy
\leq \frac{C}{n^{5 \alpha}} f(x). \nonumber
\end{eqnarray}
Denote $T = \{ x \in K_0 ; B(x, n^{-3 \alpha}) \subset K_0 \}$. We
have shown that (\ref{eq_738}) holds for any $x \in T$. Thus,
\begin{equation}
 \int_T |g(x) - f(x)| dx \leq \frac{C}{n^{5 \alpha}}
\int_T f(x) dx\leq \frac{C}{n^{5 \alpha}}. \label{eq_335}
\end{equation}
Note that $\RR^n \setminus T \subset (\RR^n \setminus K_0) \cup \{
x \in \RR^n ; \exists y \not \in K_0, |x-y| \leq n^{-3 \alpha}
\}$. Corollary \ref{cor_302} and Lemma \ref{lem_1208} show that
\begin{equation}
 \int_T f(x) dx = 1 - \int_{\RR^n \setminus T} f(x) dx \geq 1 - e^{-\alpha n / 8} - \frac{C}{n^{\alpha}}
\geq 1 - \frac{C^{\prime}}{n^{\alpha/10}}.
\label{eq_607}
\end{equation}
By (\ref{eq_335}) and (\ref{eq_607}),
\begin{equation}
 \int_{T} g(x) dx \geq \int_T f(x) dx - \int_T |g(x) - f(x)| dx
\geq  1 - \frac{\tilde{C}}{n^{\alpha / 10}}. \label{eq_608}
\end{equation} Since $\int f = \int g = 1$, then according
to (\ref{eq_607}) and (\ref{eq_608}),
\begin{equation}
 \int_{\RR^n \setminus T} |g(x) - f(x)| dx
\leq \int_{\RR^n \setminus T} \left[ g(x) + f(x) \right] dx
\leq \hat{C} n^{-\alpha
/10}.
\label{eq_336}
\end{equation}
The lemma follows by adding  inequalities
 (\ref{eq_335}) and (\ref{eq_336}).
 \hfill $\square$

\smallskip Lemma \ref{lem_915} allows us to convolve our log-concave
function with a small gaussian. The proof of the next lemma is the
most straightforward adaptation of the proof of Lemma
\ref{lem_fourier}. We  sketch the main points of difference
between the proofs.

\begin{lemma} Let $n \geq 2$ be an integer,
let $\alpha \geq 10$, and let $f: \RR^n \rightarrow [0, \infty)$
be an isotropic, log-concave function.
 Assume that
\begin{equation}
\sup_{\theta \in S^{n-1}} M_{f}(\theta, t) \leq e^{-5 \alpha n
\log n} + \inf_{\theta \in S^{n-1}} M_{f}(\theta, t) \ \ \
\text{for all} \ \ t \in \RR. \label{eq_934}
\end{equation}
Denote $g = f * \gamma_{n, n^{-\alpha}}$, where $*$ stands for
convolution. Then,
$$
\sup_{\theta \in S^{n-1}} g(t \theta) \leq e^{-\alpha n \log n} +
\inf_{\theta \in S^{n-1}} g(t \theta) \ \ \ \text{for all} \ \ t
\geq 0. $$ \label{lem_fourier2}
\end{lemma}

\emph{Sketch of proof:} For $\xi_1, \xi_2 \in \RR^n$ with $|\xi_1|
= |\xi_2| = r$,
$$
\left| \hat{f}(\xi_1) - \hat{f}( \xi_2) \right| \leq 2 \pi r
\int_{-\infty}^{\infty} \left| M_f \left( \frac{\xi_1}{|\xi_1|}, t
\right) -M_f \left( \frac{\xi_2}{|\xi_2|}, t \right) \right| dt
$$
and consequently $\left| \hat{f}(\xi_1) - \hat{f}( \xi_2) \right|
\leq r e^{-2 \alpha n \log n}$, by (\ref{eq_934}) and Lemma
\ref{lem_1028}.  Note that $\hat{g}(\xi) = \hat{f}(\xi) \cdot \exp
(-2 \pi^2 n^{-\alpha} |\xi|^2 )$ (see, e.g., \cite[page
6]{stein_weiss}). Therefore
\begin{equation}
\left| \hat{g}(\xi_1) - \hat{g}( \xi_2) \right| \leq r e^{-2 \pi^2
n^{-\alpha} r^2 } e^{-2 \alpha n \log n} \ \ \ \text{when} \ \
|\xi_1| = |\xi_2| = r. \label{eq_950}
\end{equation}
Let $x \in \RR^n$ and  $U \in O(n)$. From (\ref{eq_950}),
\begin{eqnarray}
\label{eqn_959} \nonumber \lefteqn{ \left| \int_{\RR^n} \left(
\hat{g}(\xi) - \hat{g}(U \xi) \right) e^{2 \pi i \langle x, \xi
\rangle} d\xi \right| \leq e^{-2 \alpha n \log n} \int_{\RR^n}
|\xi| e^{-2 \pi^2 n^{-\alpha} |\xi|^2 } d \xi } \\ &  =&  e^{-2
\alpha n \log n} n^{\frac{\alpha (n+1)}{ 2}} \int_{\RR^n} |\xi|
e^{-2 \pi^2 |\xi|^2 } d \xi \leq e^{-\alpha n \log n}.
\phantom{aaaaaaaaaaa}
\end{eqnarray}
Since $x \in \RR^n$ and $U \in O(n)$ are arbitrary, the lemma
follows from (\ref{eqn_959}) by the Fourier inversion formula.
\hfill $\square$

\smallskip Later, we will combine the following proposition with Lemma
\ref{lem_dvo} in order to show that a typical marginal is very
close, in the total-variation metric, to a spherically-symmetric
concentrated distribution. A random vector $X$ in $\RR^n$ has a
spherically-symmetric distribution if  $Prob \{ X \in U (A) \} =
Prob \{ X \in A \}$ for any measurable set $A \subset \RR^n$ and
an orthogonal transformation $U \in O(n)$.

\begin{proposition}
There exist universal constants $C_1,c, C > 0$ for which the
following holds: Let $n \geq 2$ be an integer, and let $f: \RR^n
\rightarrow [0, \infty)$ be an isotropic, log-concave function.
Let $X$ be a random vector in $\RR^n$ with density $f$. Assume
that
\begin{equation}
\sup_{\theta \in S^{n-1}} M_{f}(\theta, t) \leq e^{-C_1 n \log n}
+ \inf_{\theta \in S^{n-1}} M_{f}(\theta, t) \ \ \ \text{for all}
\ \ t \in \RR. \label{eq_848__}
\end{equation}
Then there exists a random vector $Y$ in $\RR^n$ such that
\begin{enumerate}
\item[(i)] $\displaystyle d_{TV}(X, Y) \leq C / n^{10} $.
\vspace{5pt} \item[(ii)]  $Y$ has a spherically-symmetric
distribution.
 \vspace{3pt} \item[(iii)] $\displaystyle
 Prob \left \{ \big| \, |Y| - \sqrt{n} \, \big| \geq \eps
 \sqrt{n}
\right \} \leq C e^{-c \eps^2 n } $  for any $0 \leq \eps \leq 1$.
\end{enumerate}
 \label{prop_842}
\end{proposition}

\emph{Proof:} Recall that
\begin{equation} Vol(\sqrt{n} D^n) \leq \hat{C}^n
\label{eq_1040}\end{equation} for some universal constant $\hat{C}
> 1$. We will define two universal constants:
$$ \alpha_0 = 10^4
[\log (\hat{C}) +1] \ \ \ \text{and} \ \ \ C_1 = \max \{ 5
\alpha_0, 2 C_0 \}
$$ where $C_0$ is the constant from Proposition
\ref{cor_735} and $\hat{C}$ is the constant from (\ref{eq_1040}).
Throughout this proof, $\alpha_0, C_0, C_1$ and $\hat{C}$ will
stand for the universal constants just mentioned. We assume that
inequality (\ref{eq_848__}) -- the main assumption of this
proposition -- holds, with the constant $C_1$ as was just defined.
We may apply Proposition \ref{cor_735}, based on (\ref{eq_848__}),
since $C_0 n \leq C_1 n \log n$. By the conclusion of that
proposition,
\begin{equation}
 Prob \left \{ \big| \ |X| - \sqrt{n} \ \big| \geq \eps \sqrt{n} \right \} \leq
C e^{-c \eps^2 n} \ \ (0 \leq \eps \leq 1). \label{eq_954}
\end{equation}  Let $Z^{\prime}$ be a gaussian random vector in $\RR^n$,
independent of $X$, with $\EE Z^{\prime} = 0$ and $Cov(Z^{\prime})
= n^{-\alpha_0} Id$. Then $\EE |Z^{\prime}|^2 = n^{1-\alpha_0}$,
and, for example, by Lemma \ref{lem_1010}(i), we know that
$$
Prob \left \{ |Z^{\prime}| \geq 1 \right \} \leq
 Prob \left \{ |Z^{\prime}| \geq 20 n \cdot \sqrt{n^{1-\alpha_0}} \right \} \leq
 e^{-n}.
 $$
Consequently, the event $-1 \leq |X + Z^{\prime}| - |X| \leq 1$
holds with probability greater than $1 - e^{-n}$. By applying
(\ref{eq_954}) we obtain that for $0 \leq \eps \leq 1$,
\begin{eqnarray}
\label{eq_1002} \lefteqn{ Prob \left \{ \big| \ |X + Z^{\prime}| -
\sqrt{n} \ \big| \geq \eps \sqrt{n} \right \} } \\ & \leq &
\nonumber e^{-n} + Prob \left \{ \big| \ |X| - \sqrt{n} \ \big|
\geq \left( \eps - \frac{1}{\sqrt{n}} \right) \sqrt{n} \right \}
\leq C^{\prime} e^{-c^{\prime} \eps^2 n}
\end{eqnarray}
(in obtaining the last inequality in (\ref{eq_1002}), one needs to
consider separately the cases $\eps < 2 / \sqrt{n}$ and $\eps \geq
2 / \sqrt{n}$).

\smallskip The density of $Z^{\prime}$ is $\gamma_{n, n^{-\alpha_0}}$.
Denote by $g = f * \gamma_{n,n^{-\alpha_0}}$ the density of the
random vector $X + Z^{\prime}$. Since $C_1 \geq 5 \alpha_0$ and
$\alpha_0 \geq 10$, then (\ref{eq_848__}) implies the main
assumption of Lemma \ref{lem_fourier2} for $\alpha = \alpha_0$. By
the conclusion of that lemma, for all $\theta_1, \theta_2 \in
S^{n-1}$ and $r \geq 0$,
\begin{equation} \left| g(r \theta_1) - g(r \theta_2) \right| \leq e^{-\alpha_0 n \log n}.
\label{eq_1122}
\end{equation}
Denote, for $x \in \RR^n$,
$$ \tilde{g}(x) = \int_{S^{n-1}} g(|x| \theta)
d\sigma_{n-1}(\theta), $$ the spherical average of $g$. The
function $\tilde{g}$ is a spherically-symmetric function with
$\int \tilde{g} = 1$, and from (\ref{eq_1122}),
\begin{equation}
 \left| \tilde{g}(x) - g(x) \right| \leq e^{-\alpha_0 n \log n} \
\ \ \text{for all} \ \ x \in \RR^n. \label{eq_1014}
\end{equation}
According to (\ref{eq_1014}) and the case $\eps = 1$ in
(\ref{eq_1002}),
\begin{eqnarray}
\label{eq_1043} \nonumber \lefteqn{ \left \| \, \tilde{g} \,  -
\, g \, \right \|_{L^1(\RR^n)} \leq \int_{|x|
\leq 2 \sqrt{n}} |\tilde{g}(x) - g(x)| dx + 2 \int_{|x| \geq 2 \sqrt{n}} g(x) dx } \\
& \leq & Vol(2 \sqrt{n} D^n) e^{-\alpha_0 n \log n} + 2 C^{\prime} e^{-c^{\prime} n}
 \leq C^{\prime \prime} e^{-c^{\prime \prime} n},
\phantom{aaaaaaaaaaaa}
\end{eqnarray}
by the definition of $\alpha_0$, where $\| F \|_{L^1(\RR^n)} =
\int_{\RR^n} |F(x)| dx$ for any measurable function $F: \RR^n
\rightarrow \RR$.

\smallskip Let $Y$ be a random variable that is distributed
according to the density $\tilde{g}$. Then $Y$ satisfies the
conclusion (ii) of the present proposition, since $\tilde{g}$ is a
radial function. Additionally, (\ref{eq_1002}) shows that $Y$
satisfies (iii), since the random variables $|Y|$ and $|X +
Z^{\prime}|$ have the same distribution. It remains to prove (i).
To that end, we employ Lemma \ref{lem_915}. The assumptions of
Lemma \ref{lem_915} are satisfied for $\alpha = \alpha_0 / 30$,
since $\alpha_0 \geq 300$. We use (\ref{eq_1043}) and the
conclusion of Lemma \ref{lem_915} to obtain
\begin{eqnarray*}
\lefteqn{  d_{TV}(X, Y)  =  \| \, f \, - \, \tilde{g} \,
 \|_{L^1(\RR^n)} \leq \| \, \tilde{g} \, - \, g \,  \|_{L^1(\RR^n)} +
 \| \, g \, - \, f \,  \|_{L^1(\RR^n)} } \\ & \leq &
C^{\prime \prime} e^{-c^{\prime \prime} n} + C n^{-\alpha_0 /300}
\leq \tilde{C} n^{-10}, \phantom{aaaaaaaaaaaaaaaaaaaaaaaa}
\end{eqnarray*}
as $\alpha_0 \geq 3000$. This completes the proof of (i). \hfill
$\square$

\begin{lemma}
Let $1 \leq k \leq n$ be integers, let $1 \leq r \leq n$, let
$\alpha, \beta > 0$ and let $X$ be a random vector in $\RR^n$ with
a spherically-symmetric distribution. Suppose
 $E \subset \RR^n$ is a $k$-dimensional
subspace. Assume that for   $0 \leq \eps \leq 1$,
\begin{equation}
 Prob \left \{  \big| \, |X| - \sqrt{n} \,
\big| \geq \eps \sqrt{n} \right \} \leq \beta e^{-\alpha \eps^2
r}. \label{eq_317_}
\end{equation}
 Then,
$$ d_{TV} \left ( \, Proj_E(X) \,  , \, Z_E \, \right ) \leq C \frac{\sqrt{k}}{\sqrt{r}} $$
where $Z_E$ is a standard gaussian random vector in $E$, and $c, C
> 0$ are constants depending only on $\alpha$ and $\beta$.
\label{computation}
\end{lemma}

\emph{Proof:} In this proof we write $c, C, C^{\prime}, \tilde{C}$
etc. to denote various positive constants depending only on
$\alpha$ and $\beta$. We may clearly assume that $n \geq 5$ and $k
\leq n-4$, as otherwise the result of the lemma is trivial with
$C \geq 2$. Let $Y$ be a random vector, independent of $X$,
 that is distributed uniformly in $S^{n-1}$.
Let $Z_E$ be a standard gaussian vector in $E$, independent of $X$
and $Y$. We will use a quantitative estimate for Maxwell's
principle by Diaconis and Freedman \cite{diaconis}. According to
their bound,
$$
 d_{TV} \left( \, Proj_E(t Y) \, , \, \frac{t}{\sqrt{n}} Z_E \,
\right) \leq 2 (k+3) / (n-k-3),
$$
for any $t \geq 0$. Since $X$ is independent of $Y$ and $Z_E$,
then also
\begin{equation}
 d_{TV}  \left( \, Proj_E(|X| Y) \, , \, \frac{|X|}{\sqrt{n}} Z_E \,
\right) \leq 2 (k+3) / (n-k-3).
\label{eq_647}
\end{equation}
 For $t \geq 0$, the density of $t Z_E$ is the function
$x \mapsto \gamma_{k, t^2}(x) \ (x \in E)$. Lemma \ref{lem_929}
implies that $d_{TV}(t Z_E, Z_E) \leq C \sqrt{k} |t^2 - 1|$, for
some universal constant $C \geq 1$. Hence,
\begin{eqnarray}
\label{eq_723}
\lefteqn{ d_{TV} \left( \, \frac{|X|}{\sqrt{n}} Z_E \, , \, Z_E \right)
\leq \EE_X \min \left \{  C \sqrt{k} \left| \frac{|X|^2}{n} - 1 \right|, 2 \right \}
}
\\ & = & \int_0^2  Prob \left \{ C \sqrt{k}  \left| \frac{|X|^2}{n} - 1 \right| \geq t \right \} dt \leq \int_0^2 C^{\prime} e^{-c^{\prime} r t^2 / k}dt \leq \tilde{C} \sqrt{\frac{k}{r}}, \nonumber
\end{eqnarray}
where we used (\ref{eq_317_}). Note that the random vectors $X$
and $|X| Y$ have the same distribution, since the distribution of
$X$ is spherically-symmetric. By combining (\ref{eq_647}) and
(\ref{eq_723}),
$$ d_{TV} \left( \, Proj_E(X) \, , \, Z_E \, \right) \leq 2
\frac{k+3}{n-k-3} + \tilde{C} \sqrt{ \frac{k}{r}} \leq \bar{C}
\sqrt{\frac{k}{r}} $$ because $r \leq n$. This completes the
proof.
 \hfill $\square$

\smallskip We are now in a position to prove Theorem \ref{thm_multi}.
Theorem \ref{thm_multi} is directly equivalent to the following
result.

\begin{theorem} Let $n \geq 1$ and $1 \leq k \leq c
\frac{\log n}{\log \log n}$ be integers, and let $X$ be a random
vector in $\RR^n$ with an isotropic, log-concave density. Then
there exists a subset $\E \subset G_{n,k}$ with $\sigma_{n,k}(\E)
\geq 1 - e^{-c n^{0.99}}$ such that for any $E \in \E$,
$$ d_{TV} \left ( \, Proj_E(X) \, , \, Z_E \, \right
) \leq   C \sqrt{k} \cdot \sqrt{\frac{\log \log n}{\log n}}, $$
where $Z_E$ is a standard gaussian random vector in $E$, and $c, C
> 0$ are universal constants. \label{thm_405}
\end{theorem}

\emph{Proof:} We use the constant $C_1$ from Proposition
\ref{prop_842}, and the constant $c$ from Lemma \ref{lem_dvo}. We
begin as in the proof of Theorem \ref{cor_202}. Denote the density
of $X$ by $f$. Set $$ \ell = \left \lfloor \frac{c}{100 C_1}
\frac{\log n}{\log \log n} \right \rfloor. $$ We may assume that
$n$ exceeds a certain universal constant, hence $\ell \geq 1$. Fix
a subspace $E \in G_{n, \ell}$ that satisfies
\begin{equation}
\sup_{\theta \in S^{n-1} \cap E} M_{f}(\theta, t) \leq e^{-C_1
\ell \log  \ell} + \inf_{\theta \in S^{n-1} \cap E} M_{f}(\theta,
t) \ \ \ \text{for all} \ \ t \in \RR. \label{eq_1032}
\end{equation}
Denote $g = \pi_E(f)$. Then $g$ is log-concave and isotropic, and
by combining (\ref{eq_1032}) with (\ref{proj_ok}) from Section
\ref{section2},
\begin{equation}
\sup_{\theta \in S^{n-1} \cap E} M_{g}(\theta, t) \leq e^{-C_1
\ell \log \ell} + \inf_{\theta \in S^{n-1} \cap E} M_{g}(\theta,
t) \ \ \ \text{for all} \ \ t \in \RR. \label{eq_1102}
\end{equation}
We invoke Proposition \ref{prop_842}, for $\ell$ and $g$, based on
(\ref{eq_1102}). Recall that $g$ is the density of $Proj_E(X)$. By
the conclusion of Proposition \ref{prop_842}, there exists a
random vector $Y$ in $E$, with a spherically-symmetric
distribution, such that
\begin{equation}
d_{TV} \left ( \, Proj_E(X) \, , \, Y \, \right ) \leq
\frac{C}{\ell^{10}} \label{eq_1104}
 \end{equation}
and
\begin{equation}
 Prob \left \{ \big| \, |Y| - \sqrt{\ell} \, \big| \geq \eps
 \sqrt{\ell}
\right \} \leq C^{\prime} e^{-c^{\prime} \eps^2 \ell } \ \ \
\text{for} \ \ 0 \leq \eps \leq 1. \label{eq_1105} \end{equation}
Fix $1 \leq k \leq \ell$, and let $F \subset E$ be a
$k$-dimensional subspace. Since the distribution of $Y$ is
spherically-symmetric, we may apply Lemma \ref{computation} for
$n=\ell$ and $r = \ell$, based on (\ref{eq_1105}). By the
conclusion of that lemma,
$$
d_{TV} \left ( \, Proj_F(Y) \, , \, Z_F \, \right )
\leq C^{\prime \prime} \frac{\sqrt{k}}{\sqrt{\ell}},
$$
where $Z_F$ is a standard gaussian random vector in $F$.
We combine the above with (\ref{eq_1104}), and obtain
\begin{equation}
d_{TV} \left ( \, Proj_F(X) \, , \, Z_F \, \right )
\leq C^{\prime \prime} \frac{\sqrt{k}}{\sqrt{\ell}}
+ \frac{C}{\ell^{10}} \leq \tilde{C}
 \frac{\sqrt{k}}{\sqrt{\ell}}.
\label{eq_1128}
 \end{equation}
 (Note that $d_{TV}(Proj_F(X), Proj_F(Y)) \leq d_{TV}(Proj_E(X),
 Y)$.)
In summary, we have proved that whenever $E$ is an
 $\ell$-dimensional subspace that satisfies (\ref{eq_1032}),
then
 all the $k$-dimensional subspaces $F \subset E$ satisfy
 (\ref{eq_1128}).

\smallskip Suppose that $E \in G_{n,\ell}$ is a random $\ell$-dimensional
subspace. We will use Lemma \ref{lem_dvo}, for $A = C_1 \log \ell$
and $\delta = 1/100$. Note that  $\ell \leq \log n$, hence $\ell
\leq c \delta A^{-1} \log n$, by the definition of $\ell$ above.
Therefore we may safely apply Lemma \ref{lem_dvo}, and conclude
that with probability greater than $1 - e^{-c n^{0.99}}$, the
subspace $E$ satisfies (\ref{eq_1032}). Therefore, with
probability greater than $1 - e^{-c n^{0.99}}$ of selecting $E$,
all $k$-dimensional subspaces
 $F \subset E$ satisfy
(\ref{eq_1128}).

\smallskip Next, we select a random subspace $F$ inside the random
subspace $E$. That is, fix $k \leq \ell - 4$, and suppose that $F
\subset E$ is a random subspace, that is
 distributed uniformly over the grassmannian
of $k$-dimensional subspaces of $E$. Since $E$ is distributed
uniformly over $G_{n,\ell}$, it follows that $F$ is distributed
uniformly over $G_{n,k}$. We thus conclude that
 $F$ -- which is a random, uniformly distributed,
$k$-dimensional subspace in $\RR^n$ --
satisfies (\ref{eq_1128}) with probability
greater than $1 - e^{-c n^{0.99}}$. Recall that $\ell
> \bar{c} (\log n) / \log \log n$ for a universal constant $\bar{c} > 0$,
and that our only assumption about $k$ was that $1 \leq k \leq
\ell$. The theorem is therefore proved. \hfill $\square$

\smallskip
\emph{Proof of Theorem \ref{thm_multi}:} Observe that
$$  \frac{1}{\sqrt{c}} \cdot \sqrt{k} \cdot \sqrt{ \frac{\log \log n}{\log n}}  \leq
\eps, $$ under the assumptions of Theorem \ref{thm_multi}. The
theorem thus follows from Theorem \ref{thm_405}, for an
appropriate choice of a universal constant $c > 0$.
 \hfill $\square$

\smallskip
\emph{Proof of Theorem \ref{thm_basic}:} Substitute $k = 1$ and $\eps =
\sqrt{\frac{\log \log n}{c \log n}}$ in Theorem \ref{thm_multi},
for $c$ being the constant from  Theorem \ref{thm_multi}. \hfill
$\square$

\smallskip
An additional notion of distance between multi-dimensional
measures is known in the literature under the name of
``$T$-distance'' (see, e.g., \cite{meckes}, \cite{NR}). For two
random vectors $X$ and $Y$ in a subspace $E \subset \RR^n$, their
$T$-distance is defined as
$$
T(X, Y) = \sup_{\theta \in S^{n-1}, t \in \RR}
 \left| \, Prob \left \{ \langle X, \theta \rangle \leq t \right \} - Prob \left \{ \langle
Y, \theta \rangle \leq t \right \} \, \right|. $$ The $T$-distance
between $X$ and $Y$
 compares only one-dimensional marginals of $X$ and $Y$,
 hence it is weaker than the total-variation distance.
The following proposition is proved
by directly
adapting the arguments of Naor and Romik \cite{NR}.

\begin{proposition}
Let $\eps > 0$, and assume that $n > \exp(C / \eps^2)$
is an integer. Suppose that $X$ is a random vector in $\RR^n$
with an isotropic, log-concave density.
Let $1 \leq k \leq c \eps^2 n$ be an
integer, and let $E \in G_{n,k}$ be a random $k$-dimensional
subspace. Then, with probability greater than $1 - e^{-c \eps^2
n}$ of choosing $E$,
$$ T \left( Proj_E(X),  Z_{E} \right) \leq \eps, $$
where $Z_E$ is a standard gaussian random vector in the subspace $E$.
Here, $c, C > 0$ are universal constants.
\label{prop_T}
\end{proposition}

\emph{Sketch of Proof:} Let $g(x) = \int_{S^{n-1}} f(|x| \theta)
d\sigma_{n-1}(\theta) \ (x \in \RR^n)$ be the spherical average of
$f$. For $0 \leq \delta \leq 1$, set $A_{\delta} = \{ x \in \RR^n
; \left| \, |x| / \sqrt{n} - 1 \, \right| \geq \delta \}$.
According to Theorem \ref{cor_202},
\begin{equation}
 \int_{A_{\delta}} g(x) dx   = \int_{A_{\delta}} f(x) dx \leq C^{\prime}
 n^{-c^{\prime}
 \delta^2} \ \ \ \text{for} \ 0 \leq \delta \leq 1.
 \label{eqn_953}
 \end{equation}
Denote $\Phi(t) = \frac{1}{\sqrt{2 \pi}} \int_{-\infty}^t
e^{-s^2/2} ds \ (t \in \RR)$ and fix $\theta_0 \in S^{n-1}$. We
apply Lemma \ref{computation} (for $r = \log n$ and $k = 1$) based
on (\ref{eqn_953}),  to obtain the inequality
\begin{equation}
\left| \int_{S^{n-1}} M_f(\theta, t) d\sigma_{n-1}(\theta) -
\Phi(t) \right| = \left| M_g(\theta_0, t) - \Phi(t) \right| \leq
\frac{C^{\prime \prime}}{\sqrt{\log n}}, \label{e_1152}
\end{equation}
valid for any $t \in \RR$. Let us fix $t \in \RR$. By Proposition
\ref{lem_bobkov}, the function $\theta \mapsto M_f(\theta, t) \
(\theta \in S^{n-1})$ is $\hat{C}$-Lipshitz. We apply Proposition
\ref{dvoretzky} for $L = \hat{C}$ and then we use (\ref{e_1152})
to conclude that with probability greater than $1 - e^{-\bar{c}
\eps^2 n}$ of selecting $E$,
\begin{equation}
\left| M_f(\theta, t) - \Phi(t) \right| \leq \eps +
\frac{C}{\sqrt{\log n}} \leq \tilde{C} \eps \ \ \ \ \ \text{for all} \ \theta \in S^{n-1} \cap E.
\label{e_101}
\end{equation}
Here we used the fact that $k \leq c \eps^2 n$. Recall that $t \in
\RR$ is arbitrary. Let $t_i = \Phi^{-1}(\eps \cdot i)$ for
$i=1,...,\lfloor 1/\eps \rfloor,$ where $\Phi^{-1}$ is the inverse
function to $\Phi$. Then, with probability greater than $1 -
e^{-c^{\prime} \eps^2 n}$ of selecting $E$, the estimate
(\ref{e_101}) holds for all $t = t_i \ (i=1,...,\lfloor 1/\eps
\rfloor)$. By using, e.g., \cite[Lemma 6]{NR} we see that with
probability greater than  $1 - e^{-c^{\prime} \eps^2 n}$ of
selecting $E$,
\begin{equation}
\left| M_f(\theta, t) - \Phi(t) \right| < \bar{C} \eps \ \ \ \ \
\forall \theta \in S^{n-1} \cap E, \ t \in \RR. \label{e_128}
\end{equation}
The proposition follows from (\ref{e_128}) and the definition of
the $T$-distance.
 \hfill
$\square$

\smallskip {\it Remark.} At first glance, the estimates
in Proposition \ref{prop_T} seem surprisingly good: Marginals of
almost-proportional dimension are allegedly close to gaussian. The
problem with Proposition \ref{prop_T} hides, first, in the
requirement that $\eps > C / \sqrt{\log n}$, and second, in the
use of the rather weak $T$-distance.

{
}

\vfill  \hfill {\small \today}

\end{document}